\documentclass[11pt]{amsart}
\usepackage{amssymb}
\usepackage[margin=1.1in]{geometry}
\usepackage{enumerate}
\usepackage{tikz}

\newtheorem{theorem}{Theorem}[section]
\newtheorem{proposition}[theorem]{Proposition}
\newtheorem{lemma}[theorem]{Lemma}
\newtheorem{corollary}[theorem]{Corollary}
\theoremstyle{definition}
\newtheorem{definition}[theorem]{Definition}
\newtheorem{example}[theorem]{Example}
\newtheorem{remark}[theorem]{Remark}
\newtheorem{problem}[theorem]{Problem}
\theoremstyle{plain}
\newtheorem{introthm}{Theorem}

\newcommand{\R}{\mathbb{R}}
\newcommand{\N}{\mathbb{N}}
\newcommand{\cl}{\operatorname{cl}}
\newcommand{\Int}{\operatorname{Int}}
\newcommand{\supp}{\operatorname{supp}}
\newcommand{\DCP}{\mathcal{DCP}}
\newcommand{\A}{\mathcal{A}}
\newcommand{\D}{\mathcal{D}}

\begin{document}

\title[Arc connected continua and the double commutant property]{Arc connected
continua and the double commutant property}

\author{A. K. Kitover}
\address{Community College of Philadelphia, 1700 Spring Garden St.,
Philadelphia, PA 19130, USA}
\email{akitover@ccp.edu}

\date{\today}


\begin{abstract}
Let $K$ be a compact metrizable space. For $f\in C(K)$ let $M_f$ be the
largest lower semicontinuous map with $k\in M_f(k)\subseteq f^{-1}(f(k))$. We
say that $K$ has the double commutant property if for every $f\in C(K)$ every
continuous function that is constant on all sets $M_f(k)$ is constant on the
fibres of $f$; this is exactly the condition that the double commutant of every
multiplication operator on $C(K)$ is the closed algebra generated by the
operator and the identity. We show that an arc connected continuum has this
property if the points having an open locally connected neighbourhood are
dense, and we extend this to some continua that are not arc connected. We then
develop new techniques that work when such points are not dense, or do not
exist at all, and prove that several well-known continua have the property:
the Cantor comb, the Cantor fan and, more generally, the cone over every
compact metrizable space, every smooth fan with a closed set of end points,
and the Cantor grid $(C\times[0,1])\cup([0,1]\times C)$. The property passes
to open images. Examples show that the hypotheses cannot be weakened in
several directions. For individual functions we show that, on a continuum all
of whose arc components are dense, every function that is constant on no arc
satisfies the double commutant condition, and that on every smooth dendroid
and every fan the functions satisfying it are dense.
\end{abstract}

\maketitle

\section{Introduction}

The initial motivation for this paper comes from operator theory. It follows from von Neumann's classical double commutant theorem (see~\cite{N}) that if $(\Omega, \Sigma, \mu)$ is a measure space, $f$ is a real-valued function from $L^\infty(\mu)$, and $M_f$ is the operator of multiplication by $f$ on $L^2(\mu)$ then the double commutant $\{M_f\}''$ of $M_f$ coincides with the closure in the strong operator topology of the algebra generated by $M_f$ and the identity operator $I$.

De Pagter and Ricker~\cite{PaR}
extended the above result to the spaces $L^p(\mu)$, $1\le p<\infty$, over Maharam measure
spaces, and showed that on $L^\infty[0,1]$ the double commutant of the operator
of multiplication by the independent variable consists of the multiplications
by Riemann integrable functions, so it is strictly larger than the
strongly closed algebra generated by the operator and the identity, which coincides with the algebra of multiplications by functions from $C[0,1]$.

In view of these results it is natural to consider the following problem. Let $K$ be a compact Hausdorff space, let $C(K)$ be the space of all real-valued continuous functions on $K$, and for $f\in C(K)$ let $\A_f$ be the closure of the algebra generated by the multiplication operator $M_f$ on $C(K)$ and the identity. What properties of $f$ and $K$ ensure that $\{M_f\}''=\A_f$?

A significant part of this problem is the following one: which properties of $K$ ensure that $\{M_f\}''=\A_f$ for every $f\in C(K)$? Clearly this property is a topological invariant of $K$, and we say that such a $K$ has the \emph{double commutant property}, or $K\in\DCP$.

We restrict our attention to the case when $K$ is a continuum (a compact metrizable connected space), and only indicate that by virtue of the results of~\cite{PaR} the maximal ideal space of the algebra $L^\infty[0,1]$ is not in $\DCP$, and that some results concerning metrizable compact spaces with isolated points are contained in~\cite[Theorem~1.15]{Ki10}.

The results obtained by the author in~\cite{Ki10} and in~\cite{Ki13} can be briefly summarized as follows.
\begin{itemize}
\item If a metrizable compact space without isolated points is in $\DCP$, then it is a continuum.
\item Every locally connected continuum is in $\DCP$.
\item There are continua that are not in $\DCP$, e.g.\ $K=\cl\{(x,\sin(1/x)):0<x\le1\}$.
\item There are arc connected but not locally connected continua, e.g.\ the broom and the bookcase, that are in $\DCP$; it is not known whether every arc connected continuum (in particular the Cantor fan) is in $\DCP$.
\item There are continua that are not arc connected but are in $\DCP$, for example $K=([-1,0]\times[-1,1])\cup\cl\{(x,\sin(1/x)):0<x\le1\}$.
\end{itemize}

In the current paper, while not arriving at a complete description of the continua in $\DCP$, we considerably strengthen the results of~\cite{Ki10} and~\cite{Ki13}. Below we state the main results; the definitions needed for them are given here and recalled where they are used.

\subsection*{The topological reformulation}
We write $f\in\D(K)$ if $\{M_f\}''=\A_f$; thus $K\in\DCP$ if and only if $\D(K)=C(K)$. The following notions translate the problem into topology.

\begin{definition}\label{d:intro-M}
Let $K$ be a compact metrizable space and $f\in C(K)$.
\begin{enumerate}[(i)]
\item A map $\Psi$ from $K$ into the nonempty subsets of $K$ is \emph{lower semicontinuous} (lsc) if $\{k:\Psi(k)\cap V\ne\emptyset\}$ is open for every open $V\subseteq K$.
\item $M_f$ denotes the largest lsc map with $k\in M_f(k)\subseteq f^{-1}(f(k))$ for all $k\in K$ (it exists; see Definition~\ref{d:M}).
\item A function $g\in C(K)$ is \emph{$f$-admissible} if it is constant on every set $M_f(k)$.
\end{enumerate}
\end{definition}

Every function of the form $h\circ f$ is $f$-admissible. The point of the definition is the converse.

\begin{introthm}[Theorem~\ref{t:char}]\label{ti:char}
Let $K$ be a compact metrizable space and $f,g\in C(K)$. Then $M_g\in\{M_f\}''$ if and only if $g$ is $f$-admissible. Consequently $f\in\D(K)$ if and only if every $f$-admissible function is constant on the fibres of $f$, and $K\in\DCP$ if and only if this holds for every $f\in C(K)$.
\end{introthm}

The argument for Theorem~\ref{ti:char} goes back to~\cite{Ki10}, and it is closely related to Ditor's theory of averaging operators and lower semicontinuous sections~\cite{Di73} (Remark~\ref{r:ditor}). Operator theory enters only through this theorem; from Section~\ref{s:lemmas} on, all arguments are topological.

\subsection*{Local structure: locally connected points and teeth}
Our sufficient conditions are expressed in terms of two kinds of local structure.

\begin{definition}\label{d:intro-L}
For a continuum $K$, $L(K)$ denotes the set of points of $K$ that have an open locally connected neighbourhood. It is an open subset of $K$.
\end{definition}

For a locally connected continuum $L(K)=K$; for the Cantor fan $L(K)=\emptyset$. Where $L(K)$ is not dense we use instead a product structure along arcs.

\begin{definition}\label{d:tooth}
A \emph{tooth chart} of a continuum $K$ is a continuous map $h:U\times[0,1]\to K$, $U$ a nonempty open subset of a compact metrizable space, together with an interval $J_h=(0,1)$ or $J_h=(0,1]$, such that $h$ maps $U\times J_h$ homeomorphically onto an open subset $O_h$ of $K$. The arcs $h(\{y\}\times J_h)$ are the \emph{open teeth}, and $h(y,0)$ is the \emph{root} of the tooth. Given a family of tooth charts, a point $v\in K$ is a \emph{star point} if it has arbitrarily small open neighbourhoods $S$ such that every $p\in S\setminus\{v\}$ is joined to $v$ by an arc $A\subseteq S$ with $A\setminus\{v\}$ contained in one open tooth.
\end{definition}

When $U$ is a clopen subset of the Cantor set and $J_h=(0,1)$, the set $O_h$ is a \emph{Cantor box} (Section~\ref{s:lemmas}), an open piece of a Cantor fence $C\times[0,1]$ in the sense of~\cite{Lip25}. The teeth of the Cantor comb are its vertical segments, rooted on the base; the teeth of a cone $\operatorname{Cone}(Y)$ are the segments $\{y\}\times(0,1]$, all rooted at the vertex, which is a star point. Star points are thus modelled on the top of a fan~\cite{Ch67}, at which the legs of the fan meet (see~\cite{BEIJKN}). Note, however, that tooth charts require open product neighbourhoods, which a general fan need not have: in the Lelek fan the end points are dense, and no nonempty open set is of the form $O_h$ (an end point of the fan lying in $O_h$ would have to be of the form $h(y,1)$ with $J_h=(0,1]$, so the nonempty open set $h(U\times(0,1))$ would contain no end points).

\subsection*{Main results}
Our most general sufficient condition is the following; it is proved as Theorem~\ref{t:teeth} together with Theorem~\ref{t:open-map}.

\begin{introthm}\label{ti:teeth}
Let $K$ be a metrizable continuum with a family of tooth charts, let $B$ be the union of the sets $O_h$, and let $\mathcal R\subseteq K$ be closed. Assume:
\begin{enumerate}[(1)]
\item $K$ is arc connected, or every arc component $P$ of $K$ is dense or satisfies $P\subseteq\cl\Int P$;
\item $K=\cl L(K)\cup\mathcal R\cup B$, and $h(U\times\{0\})\subseteq\mathcal R\cup\cl L(K)$ for every chart $h$;
\item every point of $(\mathcal R\cup\cl L(K))\setminus(L(K)\cup B)$ has arbitrarily small open neighbourhoods $U'$ such that every component of $U'\cap T$, $T$ an open tooth, has an end point in $U'\setminus B$ (an end point of $h(\{y\}\times I)$, $I$ an interval, is a point $h(y,\alpha)$ with $\alpha$ an end point of $I$);
\item every point of $\mathcal R$ lies in $\cl L(K)$ or is a star point.
\end{enumerate}
Then $K\in\DCP$, and so is every image of $K$ under an open continuous map.
\end{introthm}

Condition (3) is a mild local condition at the points outside the teeth; (2) and (4) say that every tooth is rooted either in the closure of $L(K)$ or at a star point. Theorem~\ref{ti:teeth} contains the following special cases.
\begin{itemize}
\item With no charts at all, its hypotheses say that $L(K)$ is dense and (1) holds (Theorem~\ref{t:fat}); in particular every arc connected continuum with $L(K)$ dense is in $\DCP$ (Theorem~\ref{t:arc-dense}). This contains the locally connected continua~\cite{Ki10} and the broom and the bookcase of~\cite{Ki13}.
\item The Cantor comb $M=(C\times[0,1])\cup([0,1]\times\{0\})$, in which $L(M)$ is not dense, and the cone over every nonempty compact metrizable space, in particular the Cantor fan $F$, in which $L(F)=\emptyset$, are in $\DCP$ (Theorem~\ref{t:combfan}); this answers two questions of~\cite{Ki13}. So are the Cantor ladder and the suspension of every nonempty compact metrizable space (Remark~\ref{r:combfan}).
\item Every smooth fan with a closed set of end points is an open image of the Cantor fan~\cite{CC89}, hence in $\DCP$ (Corollary~\ref{c:smooth}).
\item If $K\times Y\in\DCP$ for some nonempty compact metrizable $Y$, then $K\in\DCP$ (Corollary~\ref{c:product}).
\end{itemize}

The Cantor grid is shown to be in $\DCP$ in this paper but is not covered by Theorem~\ref{ti:teeth}, for any family of tooth charts: the points where its arcs cross lie in no open tooth and are neither in $\cl L(G)=\emptyset$ nor star points (Remark~\ref{r:teeth}). An additional argument at these points gives the following.

\begin{introthm}[Theorem~\ref{t:grid}]\label{ti:grid}
The Cantor grid $G=(C\times[0,1])\cup([0,1]\times C)$ is in $\DCP$, although $L(G)=\emptyset$.
\end{introthm}

The hypotheses of Theorem~\ref{ti:teeth} cannot simply be dropped. By Proposition~\ref{p:ray}, no compactification of a ray, or of two disjoint open rays, with a nondegenerate remainder is in $\DCP$ (for instance the $\sin(1/x)$ continuum or a ray spiralling onto a circle), although adding the disk bounded by the circle gives a continuum in $\DCP$ (Example~\ref{e:spiraldisk}). Example~\ref{e:sin-arc} shows that arc connectedness cannot be weakened to the condition that every point can be joined by an arc to $L(K)$, and Example~\ref{e:halfsquare}, in which $L(K)$ is dense, that the condition on arc components in (1) cannot be dropped. Indecomposability does not help either: by a theorem of Fukaishi and Matsuhashi~\cite[Theorem~1.1]{FM17} and Theorem~\ref{t:open-map}, there are indecomposable, even tree-like indecomposable~\cite{MOB25}, continua not in $\DCP$ (Remark~\ref{r:indec}).

\subsection*{Individual functions}
For continua not covered by Theorem~\ref{ti:teeth}, such as solenoids, the buckethandle continuum and the Lelek fan, we show that many individual functions belong to $\D(K)$.

\begin{introthm}\label{ti:functions}
Let $K$ be a metrizable continuum.
\begin{enumerate}[(i)]
\item If all arc components of $K$ are dense, then every $f\in C(K)$ that is constant on no arc belongs to $\D(K)$ (Theorem~\ref{t:nonconst}).
\item If there are retractions $r_n$ of $K$ onto subsets in $\DCP$ with $r_n\to\mathrm{id}_K$ uniformly, then $\varphi\circ r_n\in\D(K)$ for all $\varphi\in C(K)$; in particular $\D(K)$ is dense in $C(K)$ (Proposition~\ref{p:retdense}).
\end{enumerate}
\end{introthm}

Part (i) applies to arc connected continua, solenoids and the buckethandle continuum. Part (ii) applies, through the small retractions onto trees and finite fans constructed by Fugate~\cite{Fu71b,Fu71a}, to every smooth dendroid and every fan, for instance the Lelek fan (Corollaries~\ref{c:smoothdendroid} and~\ref{c:allfans}); such approximations are what Eberhart and Fugate call approximation from within~\cite{EF71}. Density of $\D(K)$ is, however, no evidence that $K\in\DCP$: the $\sin(1/x)$ continuum admits small retractions onto arcs, so $\D(K)$ is dense for it, yet it is not in $\DCP$ (Remark~\ref{r:density}).

\subsection*{Open problems and organization}
The basic question of~\cite{Ki13}, whether every arc connected continuum is in $\DCP$, remains open; so does the question whether every dendroid, or the Lelek fan, is in $\DCP$. These and related questions on retracts, continuous images and inverse limits are collected in Section~\ref{s:problems}. Section~\ref{s:char} proves Theorem~\ref{ti:char}; Section~\ref{s:lemmas} derives collapsed levels from local structure; Section~\ref{s:main} proves Theorems~\ref{ti:teeth} and~\ref{ti:grid} and contains the examples and counterexamples; Section~\ref{s:functions} deals with individual functions, including Theorem~\ref{ti:functions}, solenoids, the buckethandle continuum and fans.

\subsection*{Notation and terminology}

\emph{Spaces and functions.} Except in the first paragraphs of the
Introduction,
$K$ is a compact metrizable space with a fixed metric $d$. $C(K)$ is the Banach
space of real-valued continuous functions on $K$ with the supremum norm. For a
set $A\subseteq K$, $\cl A$ and $\Int A$ denote its closure and interior. For
a measure $\mu \in C(K)^*$, $\supp\mu$ is its support. $C$ is the ternary
Cantor set in $[0,1]$.

\emph{Operators.} $L(C(K))$ is the algebra of bounded linear operators on
$C(K)$. For $f\in C(K)$, $M_f\in L(C(K))$ is the operator of multiplication by
$f$. For a set $\mathfrak B\subseteq L(C(K))$, $\mathfrak B'$ is its commutant
and $\mathfrak B''=(\mathfrak B')'$ its double commutant. $\A_f$ is the
closure of the algebra generated by $M_f$ and the identity $I$; it is the same
in the operator norm and in the strong operator topology~\cite[Corollary~1.2]{Ki10}.
We write $f\in\D(K)$ if $\{M_f\}''=\A_f$, and $K\in\DCP$ if $\D(K)=C(K)$.

\emph{Fibres.} For $f\in C(K)$ and $t\in\R$ we write $Z_t=f^{-1}(t)$, and
$\Phi_f(k)=Z_{f(k)}$ is the fibre of $f$ through $k\in K$. When several
functions are involved, the dependence of $Z_t$ on $f$ is clear from the
context.

\emph{Lower semicontinuous maps.} A map $\Psi$ from $K$ into the nonempty
subsets of $K$ is \emph{lower semicontinuous} (lsc) if
$\{k:\Psi(k)\cap V\ne\emptyset\}$ is open for every open $V\subseteq K$;
equivalently, for every $k_0$, every $s\in\Psi(k_0)$ and every sequence
$k_n\to k_0$ there are $s_n\in\Psi(k_n)$ (for large $n$) with $s_n\to s$. For
$f\in C(K)$, $M_f$ also denotes the largest lsc map with
$k\in M_f(k)\subseteq\Phi_f(k)$ for all $k$ (Definition~\ref{d:M}). Following~\cite{Ki10}
we use the same letter as for the multiplication operator; the sets $M_f(k)$
always appear with an argument, so no confusion should arise.

\emph{Admissible functions.} A function $g\in C(K)$ is \emph{$f$-admissible} if
it is constant on every set $M_f(k)$. Every function $h\circ f$,
$h\in C(f(K))$, is $f$-admissible. By Theorem~\ref{t:char}, $g$ is
$f$-admissible if and only if $M_g\in\{M_f\}''$; hence $f\in\D(K)$ if and only
if every $f$-admissible function is constant on the fibres of $f$.

\emph{Continua.} A \emph{continuum} is a connected compact metrizable space; it
is \emph{nondegenerate} if it has more than one point. An \emph{arc} is a
homeomorphic image of $[0,1]$; its relative interior is the image of $(0,1)$.
A space is \emph{arc connected} if any two of its points are joined by an
arc, and the \emph{arc components} of a space are its maximal arc connected
subsets. $L(K)$ (Definition~\ref{d:intro-L}) is the set of points of $K$ that have an open locally
connected neighbourhood; it is an open subset of $K$.

\emph{Dendroids and fans.} A \emph{dendroid} is an arc connected hereditarily
unicoherent continuum; any two of its points $a,b$ are joined by a unique arc
$ab$. A dendroid $K$ is \emph{smooth} if there is a point $p\in K$ (an
\emph{initial point}) such that $a_n\to a$ implies $pa_n\to pa$ in the
Hausdorff metric. A \emph{fan} is a dendroid with exactly one ramification
point $v$, its \emph{vertex} (a point that is the common end point of three
arcs, pairwise disjoint except for that point); a fan is \emph{smooth} if it
is smooth with initial point $v$, and a point of a fan is an \emph{end point}
if it is an end point of every arc containing it.

Further notions (crossing points, Cantor boxes) are defined where they are
first used.

\section{Lower semicontinuous submaps and the double commutant}\label{s:char}

\begin{definition}\label{d:M}
Let $\mathfrak M_f$ be the family of lsc maps $\Psi$ with
$k\in\Psi(k)\subseteq\Phi_f(k)$ for all $k$. Put
$M_f(k)=\cl\bigcup_{\Psi\in\mathfrak M_f}\Psi(k)$.
\end{definition}

A union of lsc maps is lsc, and so is the closure of an lsc map; since the
fibres are closed, $M_f\in\mathfrak M_f$, and $M_f$ is the largest element of
$\mathfrak M_f$. Thus, to show that a point $s$ lies in $M_f(u)$, it suffices
to find some $\Psi\in\mathfrak M_f$ with $s\in\Psi(u)$.

\subsection*{The double commutant}
Let $C(K)$ be the Banach space of real-valued continuous functions on $K$. For
$f\in C(K)$ the closure of the algebra generated by $M_f$ and $I$ is the same in
the operator norm and in the strong operator topology~\cite[Corollary~1.2]{Ki10};
we denote it by $\A_f$.

\begin{lemma}[{\cite[Lemmas~1.1, 1.8]{Ki10}}]\label{l:kernel}
Every $T\in L(C(K))$ has the form $(Th)(k)=\int h\,d\mu_k$, where
$k\mapsto\mu_k\in C(K)^*$ is weak$^*$ continuous and bounded, and conversely.
Moreover $TM_f=M_fT$ if and only if $\supp|\mu_k|\subseteq\Phi_f(k)$ for all
$k\in K$, and every operator in $\{M_f\}''$ is a multiplication operator $M_g$,
$g\in C(K)$.
\end{lemma}

\begin{proof}
Only the commutation criterion needs a comment: $(TM_fh)(k)-(M_fTh)(k)=\int
(f-f(k))h\,d\mu_k$, which vanishes for all $h$ iff $(f-f(k))\mu_k=0$, i.e.\ iff
$|\mu_k|$ is concentrated on the closed set $\Phi_f(k)$.
\end{proof}

If $g\in C(K)$ is constant on the fibres of $f$, then $h(f(k)):=g(k)$ is a
well-defined function on $f(K)$ with $g=h\circ f$, and $h$ is continuous:
for a closed set $F\subseteq\mathbb R$ we have
$h^{-1}(F)=f\bigl(g^{-1}(F)\bigr)$, which is compact, being a continuous
image of the compact set $g^{-1}(F)$. (Equivalently, $f:K\to f(K)$ is a
closed continuous surjection and hence a quotient map;
see~\cite[\S2.4 and Theorem~3.1.12]{En89}.) Thus the functions in $C(K)$
constant on the fibres of $f$ are exactly the functions $h\circ f$,
$h\in C(f(K))$; by the Weierstrass
theorem
\begin{equation}\label{e:Af}
\A_f=\{M_g:\ g \text{ is constant on every fibre of } f\}.
\end{equation}

\begin{lemma}[{\cite[Proposition~3.1]{Di73}}]\label{l:supp}
Let $k\mapsto\mu_k\in C(K)^*$ be weak$^*$ continuous. Then
$k\mapsto\{k\}\cup\supp|\mu_k|$ is lsc.
\end{lemma}

\begin{proof}
Let $V$ be open and $V\cap\supp|\mu_{k_0}|\ne\emptyset$. Then $\mu_{k_0}$ does
not vanish on the functions supported in $V$, so there is $h\in C(K)$ with
$\supp h\subseteq V$ and $\mu_{k_0}(h)\ne0$. By weak$^*$ continuity
$\mu_k(h)\ne0$ for $k$ near $k_0$, hence $\supp|\mu_k|\cap V\ne\emptyset$. The
map $k\mapsto\{k\}$ is lsc, and a union of lsc maps is lsc.
\end{proof}

We shall use Michael's selection theorem in the following form.

\begin{theorem}[{Michael~\cite[Theorem~1.2]{Mi66}}]\label{t:michael}
Let $X$ be paracompact, $\Delta$ a metrizable subset of a complete locally
convex space $E$, and $\Theta$ an lsc map from $X$ into the nonempty subsets of
$\Delta$ with compact values.
Then there is a continuous $\varphi:X\to E$ with
$\varphi(x)\in\overline{\operatorname{conv}}\,\Theta(x)$ for all $x$.
\end{theorem}

We also use the following lemma due to Michael~\cite[p.~362, Examples 1.3,
1.3$^*$]{Mi56}
\begin{lemma} \label{Michael-observation}
 If $\Theta$ is lsc, $A\subseteq X$ is closed, and $s$ is a
continuous selection of $\Theta|_A$, then the map equal to $\{s(x)\}$ on $A$ and
to $\Theta$ off $A$ is lsc.
\end{lemma}

The sets $M_f(k)$ were introduced in~\cite[Definition~1.5, Lemma~1.9]{Ki10},
and the argument for the necessity part of the following theorem is contained
in the proof of~\cite[Theorem~1.14]{Ki10}; the theorem itself was not stated
there. We include the short proof for completeness.

\begin{theorem}\label{t:char}
Let $K$ be a compact metrizable space and $f,g\in C(K)$. Then
$M_g\in\{M_f\}''$ if and only if $g$ is $f$-admissible. Consequently
$f\in\D(K)$ if and only if every $f$-admissible function is constant on the
fibres of $f$.
\end{theorem}

\begin{proof}
Suppose $g$ is constant on each $M_f(k)$ and let $T\in\{M_f\}'$ have kernels
$\mu_k$. By Lemmas~\ref{l:kernel} and~\ref{l:supp},
$k\mapsto\{k\}\cup\supp|\mu_k|$ belongs to $\mathfrak M_f$, hence
$\supp|\mu_k|\subseteq M_f(k)$ and
$(TM_gh)(k)=\int gh\,d\mu_k=g(k)\int h\,d\mu_k=(M_gTh)(k)$.

Conversely, let $M_g\in\{M_f\}''$, $k_0\in K$ and $s_0\in M_f(k_0)$. Let $E$ be
the algebraic dual of $C(K)$ with the weak$^*$ topology (a complete locally
convex space) and $\Delta=\{\delta_s:s\in K\}\subseteq E$, which is homeomorphic
to $K$. The map $\Theta_0(k)=\{\delta_s:s\in M_f(k)\}$ is lsc with compact
values, and by Lemma~\ref{Michael-observation} so is $\Theta$, equal to $\{\delta_{s_0}\}$
at $k_0$ and to $\Theta_0$ elsewhere. Theorem~\ref{t:michael} gives a continuous
$\varphi:K\to E$ with $\varphi(k)$ in the closed convex hull of $\Theta(k)$.
Every such functional is positive with $\varphi(k)(\mathbf1)=1$, hence a
probability measure, and it vanishes on nonnegative functions vanishing on
$M_f(k)$; so $\supp\varphi(k)\subseteq M_f(k)\subseteq\Phi_f(k)$ and
$\varphi(k_0)=\delta_{s_0}$. The operator $(Th)(k)=\varphi(k)(h)$ is bounded on
$C(K)$ and commutes with $M_f$ by Lemma~\ref{l:kernel}; hence $TM_g=M_gT$ and
\[
g(s_0)=(TM_g\mathbf1)(k_0)=(M_gT\mathbf1)(k_0)=g(k_0).
\]
The remaining statements follow from Lemma~\ref{l:kernel}
and~\eqref{e:Af}.
\end{proof}

\begin{remark}\label{r:ditor}
Theorem~\ref{t:char} can also be obtained from Ditor's theory of regular
averaging operators and lower semicontinuous sections~\cite{Di73} (see
also~\cite{Mi64,Pe68,DE72}). Let $P_f=\{(k,s)\in K\times K:f(k)=f(s)\}$ and let
$\pi_1,\pi_2:P_f\to K$ be the coordinate projections. One checks with
Lemma~\ref{l:kernel} and~\cite[Lemma~2.1]{Di73} that the positive operators
$T\in\{M_f\}'$ with $T\mathbf1=\mathbf1$ are exactly the operators
$h\mapsto u(h\circ\pi_2)$, where $u$ is a regular averaging operator for
$\pi_1$, and that $k\mapsto\{k\}\times M_f(k)$ is the largest lsc section
$\Lambda_{\pi_1}$ of $\pi_1$ in the sense of~\cite[Section~4]{Di73}
(see~\cite[Proposition~4.2]{Di73}). The sufficiency part of
Theorem~\ref{t:char} then corresponds to~\cite[Corollary~3.2 and
Proposition~4.2(ii)]{Di73}, and the necessity part to the proof
of~\cite[Theorem~3.4]{Di73} (cf.~\cite[Remark~3.6]{Di73}), which produces a
regular averaging operator whose measures are supported in a prescribed lsc
section; the only additional point is that the section must take a prescribed
value at one point, which is provided by Michael's observation.
\end{remark}

\subsection*{Standing conventions}
From now on we use only the topological formulation. A function $f\in C(K)$
and an $f$-admissible function $g$ are fixed. A level $t$ is \emph{collapsed}
if $g$ is constant on $Z_t$, and $\Sigma$ denotes the set of points of $K$
whose level is collapsed. Thus $f\in\D(K)$ if and only if, for every
$f$-admissible $g$, every level is collapsed.

\section{Collapsed levels from local structure}\label{s:lemmas}

By Theorem~\ref{t:char}, to prove that $f\in\D(K)$ we must show that, for every
$f$-admissible $g$, every level of $f$ is collapsed. This is done in two steps.
First, in this section, we obtain a large supply of collapsed levels from local
conditions and from retractions compatible with $f$. For instance, if
$L(K)$ is dense in $K$, every level in $f(D_f)$ is collapsed for a dense set
$D_f$ (Lemma~\ref{l:dense}). Second, in Section~\ref{s:main}, we show that
under additional global assumptions on $K$, such as arc connectedness, all levels are then
collapsed; the main tool is Theorem~\ref{t:crit}, which passes the information
from collapsed levels near $c$ to the level $c$ itself.

All the conditions in the first part of this section work through the same
mechanism: they show that
$M_f(k)$ is large. Either some $k$ satisfies $M_f(k)=Z_{f(k)}$
(Lemmas~\ref{l:interior}, \ref{l:lscopen}), or a point $v$ lies in $M_f(u)$ for
\emph{every} $u$ in the fibre of $v$ (Lemma~\ref{l:open} and its special cases,
Lemmas~\ref{l:lc}, \ref{l:arc}). In both cases $g$ is constant on the fibre,
since $g$ is $f$-admissible.

\subsection*{Openness and crossing points}

\begin{lemma}[{cf.~\cite[Lemma~1.13]{Ki10}}]\label{l:interior}
If $U=\Int Z_t\ne\emptyset$, then $M_f(k)=Z_t$ for $k\in U$; in
particular the level $t$ is collapsed.
\end{lemma}

\begin{proof}
The map equal to $Z_t$ on the open set $U$ and to $\{k\}$ elsewhere belongs to
$\mathfrak M_f$.
\end{proof}

Lemma~\ref{l:interior} is the special case $U=\Int Z_t$ of the following
lemma, since $\Phi_f\equiv Z_t$ on $U$.

\begin{lemma}\label{l:lscopen}
Let $U\subseteq K$ be open and assume that $\Phi_f$ is lsc at every point of
$U$. Then $M_f(k)=\Phi_f(k)$ for every $k\in U$; in particular every
level $f(k)$, $k\in U$, is collapsed.
\end{lemma}

\begin{proof}
Let $\Psi=\Phi_f$ on $U$ and $\Psi(k)=\{k\}$ off $U$. At points of the open set
$U$, $\Psi$ coincides near the point with $\Phi_f$, hence is lsc there; at a
point $k_0\notin U$ we have $\Psi(k_0)=\{k_0\}$ and $k\in\Psi(k)$ for all $k$.
So $\Psi\in\mathfrak M_f$.
\end{proof}

The next lemma is the main tool of this section.

\begin{lemma}\label{l:open}
Let $R\subseteq K$ be such that $f|_R:R\to f(K)$ is an open map. Then $s\in M_f(u)$
for every $s\in R$ and every $u\in Z_{f(s)}$. In particular every level in
$f(R)$ is collapsed.
\end{lemma}

\begin{proof}
Let $\Psi(k)=\{k\}\cup(\Phi_f(k)\cap R)$; it suffices to show that $\Psi$ is lsc,
for then $\Psi\in\mathfrak M_f$ and $s\in\Psi(u)\subseteq M_f(u)$. Let $V$ be open
and $\Psi(k_0)\cap V\ne\emptyset$. If $k_0\in V$, then $V$ is a neighbourhood
of $k_0$, and for every $k\in V$ we have $k\in\Psi(k)\cap V$, so
$\Psi(k)\cap V\ne\emptyset$. Otherwise there is $s\in V\cap R$ with
$f(s)=f(k_0)$; since $f|_R$ is open,
$f(V\cap R)\supseteq(f(s)-\delta,f(s)+\delta)\cap f(K)$
for some $\delta>0$, and for every $k$ in the open neighbourhood
$\{k:|f(k)-f(k_0)|<\delta\}$ of $k_0$ there is $y\in V\cap R$ with $f(y)=f(k)$,
i.e.\ $y\in\Psi(k)\cap V$.
\end{proof}

\begin{remark}
With $R=K$, Lemma~\ref{l:open} gives Corollary~\ref{c:open} below in the form
$M_f(k)=\Phi_f(k)$ for all $k$. In the applications below $f|_R$ is even open
as a map into $\R$. The point of Lemma~\ref{l:open} is that $R$
may be small (for instance nowhere dense), and that nothing is assumed about
the points $u$.
\end{remark}

\begin{definition}\label{d:openness}
A point $v\in K$ is an \emph{openness point} of $f$ if $f(v)$ is an interior
point of $f(U)$ for every neighbourhood $U$ of $v$, i.e.\ $f$ is open at $v$.
\end{definition}

\begin{lemma}\label{l:lc}
Let $W\subseteq K$ be open and locally connected, and assume that $f$ is not
constant on any nonempty open subset of $W$. If $v\in W$ is an openness point of
$f$, then $v\in M_f(u)$ for every $u\in Z_{f(v)}$.
\end{lemma}

\begin{proof}
Let $O$ be the set of openness points of $f$ in $W$; by Lemma~\ref{l:open} it
suffices to show that $f|_O$ is open. Let $s\in O$ and let $V$ be an open
neighbourhood of $s$. Choose a connected open $V'\subseteq V\cap W$ containing
$s$, and $\delta>0$ with $(f(s)-\delta,f(s)+\delta)\subseteq f(V')$. Let
$|t-f(s)|<\delta$. The sets $A=V'\cap\{f<t\}$ and $B=V'\cap\{f>t\}$ are nonempty
and open, and $V'\setminus(A\cup B)=V'\cap Z_t$ has empty interior; hence $V'$ is
the union of the relative closures of $A$ and $B$, which therefore meet at a
point $y\in V'$. Then $f(y)=t$, and since $K$ is locally connected at $y$, every
neighbourhood of $y$ contains a connected neighbourhood meeting $A$ and $B$, so
$y$ is an openness point. Thus $y\in O\cap V$ and $f(y)=t$, so
$(f(s)-\delta,f(s)+\delta)\subseteq f(V\cap O)$.
\end{proof}

Lemma~\ref{l:lc} strengthens~\cite[Lemma~1.11]{Ki10} and Lemma 4
of~\cite{Ki13}: no assumption is made on $u$, and local connectedness is used
only near $v$. The same argument works along arcs, without any local
connectedness.

\begin{definition}
Let $\gamma\subseteq K$ be an arc such that $f|_\gamma$ is not constant on any
nondegenerate subarc. A point $v$ of the relative interior $\gamma^\circ$ is a
\emph{crossing point of $f$ on $\gamma$} if every subarc of $\gamma$ containing
$v$ in its relative interior contains points with $f<f(v)$ and points with
$f>f(v)$.
\end{definition}

\begin{lemma}\label{l:arc}
If $v$ is a crossing point of $f$ on an arc $\gamma$, then $v\in M_f(u)$ for
every $u\in Z_{f(v)}$.
\end{lemma}

\begin{proof}
Let $O_\gamma$ be the set of crossing points on $\gamma$; by
Lemma~\ref{l:open} it suffices to show that $f|_{O_\gamma}$ is open. Let $s\in
O_\gamma$ and let $V$ be an open neighbourhood of $s$. Choose a subarc
$J\subseteq V\cap\gamma^\circ$ with $s$ in its relative interior; then $f(J)$ is
an interval containing $f(s)$ in its interior, say
$(f(s)-\delta,f(s)+\delta)\subseteq f(J)$. For $|t-f(s)|<\delta$ the relatively
open sets $J\cap\{f<t\}$ and $J\cap\{f>t\}$ are nonempty, and $J\cap Z_t$ has
empty interior in $J$; since $J$ is connected, their relative closures meet at
some $y\in J$. Then $f(y)=t$, and every subarc of $\gamma$ containing $y$ in its
relative interior contains points of $J$ near $y$ with $f<t$ and with $f>t$; so
$y\in O_\gamma\cap V$.
\end{proof}

The next result was stated in~\cite[Theorem~1.19(1)]{Ki10}, where its proof
was omitted as being very similar to that of~\cite[Theorem~1.14]{Ki10}; it is
now an immediate consequence of Lemma~\ref{l:lscopen}.

\begin{corollary}[{\cite[Theorem~1.19(1)]{Ki10}}]\label{c:open}
If $f:K\to f(K)$ is an open map, then $f\in\D(K)$.
\end{corollary}

\begin{proof}
If $f$ is open, $t\mapsto f^{-1}(t)$ is lsc on $f(K)$ (equivalently, the
decomposition of $K$ into the fibres of $f$ is
continuous~\cite[Theorem~1.2.25]{Ma05}), so $\Phi_f$ is lsc on $K$; by Lemma~\ref{l:lscopen} with $U=K$, $M_f=\Phi_f$, so every
$f$-admissible function is constant on the fibres of $f$.
\end{proof}

Let $D_f$ be the set of points $x\in K$ such that $x\in\Int Z_{f(x)}$, or
$\Phi_f$ is lsc at every point of some neighbourhood of $x$, or $x\in R$ for some
$R\subseteq K$ with $f|_R$ open (this includes the points satisfying the
hypotheses of Lemma~\ref{l:lc} and the crossing points on arcs). By
Lemmas~\ref{l:interior}, \ref{l:lscopen} and~\ref{l:open}, the level $f(x)$ is
collapsed for every $x\in D_f$.

\begin{lemma}\label{l:dense}
If $L(K)$ is dense in $K$, then $D_f$ is dense in $K$ for every $f\in C(K)$.
\end{lemma}

\begin{proof}
Let $U$ be nonempty and open. There is a nonempty connected open locally
connected $V\subseteq U$. If $f$ is constant on a nonempty open subset of $V$,
that subset lies in $D_f$. Otherwise $f(V)$ is a nondegenerate interval; for
$t$ in its interior the argument in the proof of Lemma~\ref{l:lc} gives an
openness point $y\in V$, and $y$ satisfies the hypotheses of Lemma~\ref{l:lc}
with $W=V$.
\end{proof}

\subsection*{Local retractions and retractions}

\begin{lemma}\label{l:locret}
Let $O\subseteq K$ be open and $r:O\to K$ continuous with $f\circ r=f$ on $O$.
Then $r(k)\in M_f(k)$ for every $k\in O$. In particular, $g\circ r=g$ on $O$.
\end{lemma}

\begin{proof}
The map equal to $\{k,r(k)\}$ on $O$ and to $\{k\}$ off $O$ is lsc and belongs
to $\mathfrak M_f$.
\end{proof}

\begin{lemma}\label{l:retract}
Let $L$ be a closed subset of $K$, let $r:K\to L$ be a retraction, and let
$f=f\circ r$. Then $r(k)\in M_f(k)$ for every $k\in K$ and
$M_{f|_L}(l)\subseteq M_f(l)$ for every $l\in L$. Consequently $g=g\circ r$,
and $g|_L$ is $f|_L$-admissible.
\end{lemma}

\begin{proof}
The first statement is Lemma~\ref{l:locret} with $O=K$. Put
$\Psi(k)=\{k\}\cup M_{f|_L}(r(k))$. Since $f(M_{f|_L}(r(k)))=\{f(r(k))\}=
\{f(k)\}$, we have $\Psi(k)\subseteq\Phi_f(k)$. For $V$ open in $K$,
$\{k:M_{f|_L}(r(k))\cap V\ne\emptyset\}=r^{-1}(\{l\in L:M_{f|_L}(l)\cap(V\cap L)
\ne\emptyset\})$ is open, because $M_{f|_L}$ is lsc on $L$ and $r$ is
continuous; hence $\Psi$ is lsc and $\Psi\in\mathfrak M_f$. For $l\in L$ we get
$M_{f|_L}(l)\subseteq\Psi(l)\subseteq M_f(l)$.
\end{proof}

\begin{corollary}\label{c:retract}
Let $r:K\to L$ be a retraction and $f=f\circ r$. If $f|_L\in\D(L)$, in
particular if $L\in\DCP$, then $f\in\D(K)$. This applies whenever $L$ is a
locally connected continuum, e.g.\ an arc or a dendrite  (a locally connected
continuum containing no simple closed curve).
\end{corollary}

\begin{proof}
By Lemma~\ref{l:retract} and the assumption, $g|_L=h\circ f|_L$ for some
$h\in C(f(L))$. Since $f=f\circ r$, we have $f(L)=f(K)$ and
$g=g\circ r=h\circ f\circ r=h\circ f$. Locally connected continua are in
$\DCP$ by~\cite[Theorem~1.14]{Ki10}; see also Corollary~\ref{c:old}.
\end{proof}

The converse of Lemma~\ref{l:retract} fails: a function $g$ with $g=g\circ r$
and $g|_L$ $f|_L$-admissible need not be $f$-admissible
(Example~\ref{e:conv}). In particular, the lemma does not show that $\DCP$
passes from $K$ to its retracts.

\subsection*{Cantor boxes}

A \emph{Cantor box} in $K$ is an open set $B$ with a homeomorphism
$B\cong A\times J$, $A$ a nonempty clopen subset of the Cantor set $C$ and $J$
an open interval; $f_x(\tau)=f(x,\tau)$. A point $(x,\tau)\in B$ is a
\emph{crossing point} if every subinterval of $J$ containing $\tau$ in its
interior contains points where $f_x<f(x,\tau)$ and points where
$f_x>f(x,\tau)$. (Unlike the definition for arcs, $f_x$ may be constant on
nondegenerate subintervals of $J$; only the behaviour of $f_x$ near $\tau$
matters.) This property depends only on $f$ along the arc of $B$ through
$(x,\tau)$ near this point, so it does not depend on the box; in particular a
crossing point of $B$ is a crossing point of every smaller box
$A'\times J'\subseteq B$ containing it. The two notions of crossing point agree
when both apply: if $\gamma$ is an arc such that $f|_\gamma$ is not constant on
any nondegenerate subarc and $v\in\gamma^\circ\cap B$, then near $v$ the arc
$\gamma$ lies in the arc of $B$ through $v$, so $v$ is a crossing point of $f$
on $\gamma$ if and only if it is a crossing point of $B$.

\begin{lemma}\label{l:baire}
Let $B=A\times J$ be a Cantor box and $t\in\R$ such that for every $x\in A$,
$f_x$ takes values $<t$ and $>t$ on $J$. Then either some $(x,\tau)\in B$ with
$f(x,\tau)=t$ is a crossing point, or $Z_t\cap B$ has nonempty interior.
\end{lemma}

\begin{proof}
If $x$ has no crossing point at level $t$, pick $a,b\in J$ with $f_x(a)<t<f_x(b)$,
say $a<b$; let $\sigma=\inf\{\tau\in[a,b]:f_x(\tau)>t\}$ and
$\rho=\sup\{\tau\in[a,\sigma]:f_x(\tau)<t\}$. If $\rho=\sigma$ this point is a
crossing point, so $\rho<\sigma$ and $f_x\equiv t$ on $[\rho,\sigma]$. Hence
$A=\bigcup E_{p,q}$ over rationals $p<q$ in $J$, where
$E_{p,q}=\{x\in A:f_x\equiv t\text{ on }[p,q]\}$ is closed. By Baire's theorem
some $E_{p,q}$ contains a clopen $A'\ne\emptyset$, and $A'\times(p,q)\subseteq
Z_t\cap B$ is open.
\end{proof}

\begin{proposition}\label{p:box}
Let $W\subseteq K$ be open, covered by Cantor boxes, and assume $f$ is not
constant on any nonempty open subset of $W$. Let $O_W$ be the set of crossing
points in $W$. Then $f|_{O_W}$ is open; hence every level in $f(O_W)$ is
collapsed (Lemma~\ref{l:open}).
\end{proposition}

\begin{proof}
Let $s=(x_0,\tau_0)\in O_W$ and $V$ a neighbourhood of $s$. Choose a box
$A'\times J'\subseteq V\cap W$ about $s$ and values $f_{x_0}(\alpha)<f(s)<
f_{x_0}(\beta)$ with $\alpha,\beta\in J'$; shrinking $A'$, for some $\delta>0$
every $x\in A'$ and $|t-f(s)|<\delta$ satisfy $f_x(\alpha)<t<f_x(\beta)$. By
Lemma~\ref{l:baire} and the hypothesis on $W$, there is a crossing point $y\in
A'\times J'$ with $f(y)=t$.
\end{proof}

\begin{remark}\label{r:box}
Zero-dimensionality of $A$ is used in Lemma~\ref{l:baire} and
Proposition~\ref{p:box} only through the existence of small clopen sets. Call
a \emph{box} an open set $B\subseteq K$ with a homeomorphism $B\cong U\times
J$, $U$ a nonempty open subset of a compact metrizable space and $J$ an open
interval, and define its crossing points as above, along the arcs
$\{y\}\times J$. Then Lemma~\ref{l:baire} and Proposition~\ref{p:box} hold for
boxes, with ``clopen'' replaced by ``open'': $U$ is a Baire space, so some
closed set $E_{p,q}$ contains a nonempty open $U''\subseteq U$, and
$U''\times(p,q)$ is open in $K$. Since for $U$ not zero-dimensional the arcs
$\{y\}\times J$ depend on the product structure, in Proposition~\ref{p:box}
one takes $O_W$ to be the set of points that are crossing points of some box
contained in $W$ (each with its given product structure).
\end{remark}

\section{Continua in $\DCP$: sufficient conditions and examples}\label{s:main}

In this section we pass from the local information of
Section~\ref{s:lemmas} to global conclusions about $K$. Theorem~\ref{t:crit}
gives sufficient conditions for a single level $c$ to be collapsed: (a) the
values of $g$ on $Z_c$ can be approached through collapsed levels above,
resp.\ below, $c$; (b) a topological condition on how $Z_c$ is reached from
$\{f>c\}$ and $\{f<c\}$, which holds automatically when $K$ is arc connected.
Combined with Lemma~\ref{l:dense}, this shows that every arc connected
continuum $K$ with $L(K)$ dense is in $\DCP$ (Theorem~\ref{t:arc-dense}); we
then replace arc connectedness by a condition on arc components
(Theorem~\ref{t:fat}). After examples, which include the earlier results
of~\cite{Ki10,Ki13} (Corollary~\ref{c:old}), we show that $\DCP$ passes to open
images (Theorem~\ref{t:open-map}) and that the Cantor comb, all cones over
compact metrizable spaces (in particular the Cantor fan) and the Cantor grid
are in $\DCP$, although $L(K)$ need not be dense for them; in
the fan and the grid no point has an open locally connected neighbourhood
(Theorems~\ref{t:combfan} and~\ref{t:grid}); Theorem~\ref{t:teeth} puts most
of these results into one framework. The section ends with
counterexamples.

We first fix some notation. Recall that $\Sigma$ is the set of points of $K$
whose level is collapsed. For $c\in f(K)$ put $Z_c^+=Z_c\cap\cl\{f>c\}$ and $Z_c^-=Z_c\cap\cl\{f<c\}$, and
\[
Z_c^{\uparrow}=\{a\in Z_c:\ \forall\ \text{neighbourhood } U \text{ of } a\
\exists\,\varepsilon>0:\ [c,c+\varepsilon)\subseteq f(U)\},
\]
and $Z_c^{\downarrow}$ similarly with $(c-\varepsilon,c]$. Clearly
$Z_c^\uparrow\subseteq Z_c^+$ and $Z_c^\downarrow\subseteq Z_c^-$; if
$\Int Z_c=\emptyset$ then $Z_c=Z_c^+\cup Z_c^-$.

\begin{lemma}\label{l:match}
Let $a\in Z_c^\uparrow$ and $u\in Z_c$, and suppose that there are
$y_n\in\cl\Sigma$ with $f(y_n)>c$, $f(y_n)\to c$ and $g(y_n)\to g(u)$. Then
$g(u)=g(a)$. The same holds with $Z_c^\downarrow$ and $f(y_n)<c$.
\end{lemma}

\begin{proof}
Choose $q_n\in\Sigma$ so close to $y_n$ that $f(q_n)>c$,
$|f(q_n)-f(y_n)|<1/n$ and $|g(q_n)-g(y_n)|<1/n$; then $f(q_n)\to c$ and
$g(q_n)\to g(u)$. Let $U_m$ be the ball of radius $1/m$ about $a$ and
$\varepsilon_m$ as in the definition of $Z_c^\uparrow$. Choose $n_1<n_2<\dots$
with $f(q_{n_m})<c+\varepsilon_m$ and $a_m\in U_m$ with $f(a_m)=f(q_{n_m})$.
Since this level is collapsed, $g(a_m)=g(q_{n_m})$; letting $m\to\infty$ gives
$g(a)=g(u)$.
\end{proof}

\begin{theorem}\label{t:crit}
Let $K$ be a connected compact metrizable space, $f\in C(K)$, $g$
$f$-admissible, and let $c\in f(K)$ with $\Int Z_c=\emptyset$. Assume:
\begin{enumerate}[(a)]
\item for every $u\in Z_c^+$ (resp.\ $u\in Z_c^-$) there are $y_n\in\cl\Sigma$
with $f(y_n)>c$ (resp.\ $f(y_n)<c$), $f(y_n)\to c$ and $g(y_n)\to g(u)$;
\item $Z_c\subseteq Y^+\cup Y^-$, where $Y^+=Z_c^+$ if $Z_c^\uparrow\ne\emptyset$
and $Y^+=\emptyset$ otherwise, and $Y^-$ is defined similarly.
\end{enumerate}
Then $g$ is constant on $Z_c$.
\end{theorem}

\begin{proof}
If $Y^+\ne\emptyset$, fix $a\in Z_c^\uparrow$; by (a) and Lemma~\ref{l:match},
$g\equiv g(a)$ on $Z_c^+$. Similarly $g\equiv g(b)$ on $Z_c^-$ for
$b\in Z_c^\downarrow$ if $Y^-\ne\emptyset$. If only one of $Y^\pm$ is nonempty
we are done by (b). If both are nonempty, then $\{f>c\}$ and $\{f<c\}$ are
nonempty; since $K\setminus Z_c$ is dense, $K=\cl\{f>c\}\cup\cl\{f<c\}$, and by
connectedness these closed sets meet, necessarily in $Z_c^+\cap Z_c^-$. Hence
$g(a)=g(b)$.
\end{proof}

Condition (a) holds in particular if every point of $Z_c^+$ is the limit of
points $u_n\in\Sigma$ with $f(u_n)>c$, and similarly for $Z_c^-$ (take
$y_n=u_n$). This is how (a) is verified in all our results except
Theorems~\ref{t:combfan}, \ref{t:teeth} and~\ref{t:grid} on the Cantor comb
and fan, rooted teeth and the Cantor grid, where the more general form is
needed (Lemmas~\ref{l:tooth-approx} and~\ref{l:approx}).

If $K$ is arc connected, condition (b) holds for every $f$ and $c$: if
$Z_c^+\ne\emptyset$, take an arc from a point with $f>c$ to a point of $Z_c$ and
let $a$ be its first point in $Z_c$; then $f>c$ on the part of the arc before
$a$, so $a\in Z_c^\uparrow$. Hence $Y^\pm=Z_c^\pm$ and $Z_c=Y^+\cup Y^-$.

\begin{remark}\label{r:updown}
Example~\ref{e:sin-arc} below explains the role of the sets
$Z_c^{\uparrow},Z_c^{\downarrow}$ in (b): there, at the level $c=0$ the fibre
$E$ meets $\cl\{f>0\}$, but no point of $E$ is approached from above along a
connected set, so $Z_0^\uparrow=\emptyset$, condition (b) fails, and indeed
$f\notin\D(K)$.
\end{remark}

\begin{theorem}\label{t:arc-dense}
Let $K$ be an arc connected metrizable continuum such that the set $L(K)$ of
points having an open locally connected neighbourhood is dense in $K$. Then
$K\in\DCP$.
\end{theorem}

\begin{proof}
Let $f\in C(K)$, let $g$ be $f$-admissible and $c\in f(K)$. If
$\Int Z_c\ne\emptyset$, use Lemma~\ref{l:interior}. Otherwise (b) holds by the
remark above, and (a) holds because $D_f\subseteq\Sigma$ is dense
(Lemma~\ref{l:dense}) and $\{f>c\}$, $\{f<c\}$ are open. By
Theorem~\ref{t:crit} every level is collapsed.
\end{proof}

\begin{corollary}\label{c:exhaust}
Let $K$ be a metrizable continuum and $K_1\subseteq K_2\subseteq\cdots$ locally
connected subcontinua of $K$ such that $K_m\subseteq\cl\Int K_m$ for every $m$,
$\bigcup_mK_m$ is dense in $K$, and every point of $K$ is joined by a path in
$K$ to a point of $\bigcup_mK_m$. Then $K\in\DCP$.
\end{corollary}

\begin{proof}
The set $\Int K_m$ is open in $K$ and is an open subset of the locally
connected space $K_m$, hence locally connected; so $\Int K_m\subseteq L(K)$.
Since $\Int K_m$ is dense in $K_m$ and $\bigcup_mK_m$ is dense in $K$, $L(K)$
is dense in $K$. Each $K_m$ is arc connected~\cite[Theorem~5.1, p.~36]{Why},
hence so is the increasing union $\bigcup_mK_m$; since every point of $K$ is
joined to it by a path, hence by an arc, $K$ is arc connected. Apply
Theorem~\ref{t:arc-dense}.
\end{proof}

\subsection*{Arc components}

\begin{lemma}\label{l:dense-arc}
Let $K$ be a continuum such that every arc component $P$ of $K$ is either dense
in $K$ or satisfies $P\subseteq\cl\Int P$. Then condition (b) of
Theorem~\ref{t:crit} holds for every $f\in C(K)$ and every $c\in f(K)$ with
$\Int Z_c=\emptyset$.
\end{lemma}

\begin{proof}
We first show: \emph{if $Z_c^+\ne\emptyset$ and $Z_c^\uparrow=\emptyset$, then
$Z_c^\downarrow\ne\emptyset$ and $Z_c^+\subseteq Z_c^-$.}

Let $W=\{f>c\}$. If an arc starting at a point of $W$ leaves $W$, its first
point outside $W$ lies in $Z_c$ and is approached from above along the arc,
hence lies in $Z_c^\uparrow$. Therefore $W$ is a union of arc components of
$K$. Let $u\in Z_c^+$; then $u\in\cl W\setminus W$, so the arc component $P_u$
of $u$ is disjoint from $W$. If $P_u$ were dense, the closed set $K\setminus
W\supseteq P_u$ would be $K$, contradicting $W\ne\emptyset$. Hence
$u\in\cl\Int P_u$. Since $f\le c$ on the open set $\Int P_u$ and
$\Int Z_c=\emptyset$, there are points $v\in\Int P_u$ arbitrarily close to $u$
with $f(v)<c$; thus $u\in Z_c^-$. Moreover, on an arc in $P_u$ from such a $v$
to $u$ the first point of $Z_c$ lies in $Z_c^\downarrow$.

By symmetry, if $Z_c^-\ne\emptyset$ and $Z_c^\downarrow=\emptyset$, then
$Z_c^\uparrow\ne\emptyset$ and $Z_c^-\subseteq Z_c^+$. Now recall
$Z_c=Z_c^+\cup Z_c^-$. If $Z_c^\uparrow$ and $Z_c^\downarrow$ are both
nonempty, then $Y^\pm=Z_c^\pm$ and (b) holds. If $Z_c^\uparrow=\emptyset$,
then either $Z_c^+=\emptyset$, and $Z_c=Z_c^-=Y^-$ because $Z_c^\downarrow\ne
\emptyset$ by the symmetric statement; or $Z_c^+\ne\emptyset$, and
$Z_c\subseteq Z_c^-=Y^-$ by the statement proved above. The case
$Z_c^\downarrow=\emptyset$ is symmetric.
\end{proof}

\begin{theorem}\label{t:fat}
Let $K$ be a metrizable continuum such that $L(K)$ is dense in $K$ and every
arc component $P$ of $K$ is either dense in $K$ or satisfies
$P\subseteq\cl\Int P$. Then $K\in\DCP$.
\end{theorem}

\begin{proof}
Let $f\in C(K)$, let $g$ be $f$-admissible and $c\in f(K)$. If
$\Int Z_c\ne\emptyset$, use
Lemma~\ref{l:interior}. Otherwise (a) holds as in the proof of
Theorem~\ref{t:arc-dense} (Lemma~\ref{l:dense}), and (b) holds by
Lemma~\ref{l:dense-arc}. Apply Theorem~\ref{t:crit}.
\end{proof}

\subsection*{Examples}

\begin{corollary}\label{c:old}
The following continua are in $\DCP$: locally connected metrizable continua
\cite[Theorem~1.14]{Ki10}; the broom and the bookcase~\cite{Ki13}; the continua
of~\cite[Proposition~7]{Ki13}.
\end{corollary}

\begin{proof}
Locally connected continua are arc connected~\cite[Theorem~5.1, p.~36]{Why}.
For the broom $L(K)$ contains the segments $\{y=x/n,\ 0<x\}$ without the
origin, for the bookcase it contains the
open horizontal segments at the heights $1/n$ (but not the bottom segment
$y=0$, which is a limit of them) and the vertical sides minus the bottom
corners; both
are dense and both continua are arc connected. The continua
of~\cite[Proposition~7]{Ki13} satisfy the hypotheses of
Corollary~\ref{c:exhaust}.
\end{proof}

\begin{remark}\label{r:prop7}
Corollary~\ref{c:exhaust} is~\cite[Proposition~7]{Ki13} with fewer
hypotheses: condition~(4) there (every point of $K_m$ is a point of local
connectedness of $K$) is not needed, and of condition~(1)
($K_m\subsetneqq K_{m+1}$) only the inclusion $K_m\subseteq K_{m+1}$ is used.
The path condition~(6) cannot be omitted: for the closure $S$ of
$\{(x,\sin(1/x)):0<x\le1\}$ and $K_m$ the part of the curve over $[1/m,1]$, all
other conditions hold, but $S\notin\DCP$ (Example~\ref{e:sin}). On the other
hand, conditions (1), (2) and (6) alone imply that $K$ is arc connected. Hence
an example showing that the density conditions (3) or (5) cannot be omitted
would be an arc connected continuum outside $\DCP$, and would answer
Problem~\ref{p:arc} in the negative.
\end{remark}

Theorem~\ref{t:fat} contains Theorem~\ref{t:arc-dense} (one arc component,
$P=K$) and the following result of~\cite{Ki13}.

\begin{proposition}[{\cite[Proposition~9]{Ki13}}]\label{p:square}
Let $Q=[-1,0]\times[-1,1]$, $C_2=\{(x,\sin(1/x)):0<x\le1\}$ and $K=Q\cup C_2$.
Then $K\in\DCP$, although $K$ is not arc connected.
\end{proposition}

\begin{proof}
The arc components of $K$ are $Q$ and $C_2$; we have $Q=\cl\Int Q$, and $C_2$
is open in $K$. Moreover $L(K)\supseteq(Q\setminus E)\cup C_2$, where
$E=\{0\}\times[-1,1]$, which is dense. Apply Theorem~\ref{t:fat}.
\end{proof}

\begin{example}\label{e:cone}
Let $S$ be the closure of $\{(x,\sin(1/x)):0<x\le1\}$. The cone
$K=\operatorname{Cone}(S)$ is
arc connected (every point is joined to the vertex by a segment), and the
image of $\{(x,\sin(1/x)):0<x\le1\}\times(0,1)$ is an open locally connected
dense subset of $K$. By Theorem~\ref{t:arc-dense}, $K\in\DCP$, although
$S\notin\DCP$. (By Theorem~\ref{t:combfan} below, the cone over every compact
metrizable space is in $\DCP$.)
\end{example}

The next two continua from~\cite{Il25} are not arc connected.

\begin{example}\label{e:spiraldisk}
Let $R=\{(1+\frac1{1+t})(\cos t,\sin t):t\ge0\}$ be a ray spiralling onto the
unit circle $S^1$, let $D$ be the closed unit disk and $K=D\cup R$ (the
``spiral with a disk'' of~\cite[Exercise~13.14]{Il25}). The arc components of
$K$ are $R$, which is open in $K$, and $D=\cl\Int D$, and $L(K)\supseteq
R\cup\Int D$ is dense. By Theorem~\ref{t:fat}, $K\in\DCP$. Without the disk,
$S^1\cup R\notin\DCP$ by Proposition~\ref{p:ray} below; so adding a disk can
produce a continuum in $\DCP$.
\end{example}

\begin{example}\label{e:tornado}
\emph{Bing's double tornado}~\cite[Exercise~13.28]{Il25} is the union of two
disjoint subsets $C_1,C_2$ of $\R^3$, each homeomorphic to
$S^1\times[0,\infty)$, such that each $C_i$ spirals onto the base circle (the
image of $S^1\times\{0\}$) of the other: $\cl C_i$ is the union of $C_i$ and
the base circle of $C_{3-i}$. As in the $\sin(1/x)$ continuum, no arc joins
$C_1$ and $C_2$, so these are the arc components. The interior of $C_i$ is
$C_i$ without its base circle, so $C_i\subseteq\cl\Int C_i$, and $L(K)\supseteq
\Int C_1\cup\Int C_2$ is dense. By Theorem~\ref{t:fat}, $K\in\DCP$.
\end{example}

\subsection*{Open images}

\begin{lemma}\label{l:openmap}
Let $p:K\to D$ be an open continuous surjection of compact metrizable spaces,
$h\in C(D)$ and $f=h\circ p$. Then $M_f(k)=p^{-1}(M_h(p(k)))$ for every
$k\in K$.
\end{lemma}

\begin{proof}
Put $\Psi(k)=p^{-1}(M_h(p(k)))$. Then $k\in\Psi(k)\subseteq\Phi_f(k)$, and for
$V$ open in $K$,
$\{k:\Psi(k)\cap V\ne\emptyset\}=p^{-1}(\{y:M_h(y)\cap p(V)\ne\emptyset\})$ is
open because $p(V)$ is open and $M_h$ is lsc. Hence $\Psi\in\mathfrak M_f$ and
$\Psi(k)\subseteq M_f(k)$.

Conversely, put $\Psi'(y)=\cl\bigcup_{x\in p^{-1}(y)}p(M_f(x))$. Then
$y\in\Psi'(y)\subseteq\Phi_h(y)$. If $W\subseteq D$ is open and
$\Psi'(y_0)\cap W\ne\emptyset$, there are $x_0\in p^{-1}(y_0)$ and $w\in
M_f(x_0)$ with $p(w)\in W$. Since $M_f$ is lsc, $M_f(x)\cap p^{-1}(W)\ne
\emptyset$ for all $x$ in a neighbourhood $N$ of $x_0$, and since $p$ is open,
$p(N)$ is a neighbourhood of $y_0$; thus $\Psi'(y)\cap W\ne\emptyset$ for $y$
near $y_0$. Hence $\Psi'\in\mathfrak M_h$, so $p(M_f(k))\subseteq M_h(p(k))$
and $M_f(k)\subseteq p^{-1}(M_h(p(k)))$.
\end{proof}

\begin{theorem}\label{t:open-map}
Let $p:K\to D$ be an open continuous surjection of compact metrizable spaces,
$h\in C(D)$ and $f=h\circ p$. Then the $f$-admissible functions are exactly the
functions $g_0\circ p$ with $g_0$ $h$-admissible, and $f\in\D(K)$ if and only
if $h\in\D(D)$. In particular, if $K\in\DCP$ then $D\in\DCP$.
\end{theorem}

\begin{proof}
By Lemma~\ref{l:openmap}, $M_f(k)\supseteq p^{-1}(p(k))$, so an
$f$-admissible $g$ is constant on the fibres of $p$ and, $p$ being a closed
surjection, $g=g_0\circ p$ with $g_0\in C(D)$; moreover $g_0$ is constant on
$p(M_f(k))=M_h(p(k))$. Conversely, if $g_0$ is $h$-admissible, then $g_0\circ
p$ is constant on $p^{-1}(M_h(p(k)))=M_f(k)$. If $h\in\D(D)$, then $g_0=H\circ
h$ and $g=H\circ f$. If $f\in\D(K)$ and $g_0$ is $h$-admissible, then $g_0\circ
p=H\circ f=H\circ h\circ p$, and $g_0=H\circ h$ since $p$ is surjective.
\end{proof}

\begin{corollary}\label{c:product}
If $Y$ is a nonempty compact metrizable space and $K\times Y\in\DCP$, then
$K\in\DCP$.
\end{corollary}

\begin{proof}
The projection $K\times Y\to K$ is an open continuous surjection; apply
Theorem~\ref{t:open-map}.
\end{proof}

\begin{remark}\label{r:Ldense}
Open maps preserve the density of $L$: if $p:K\to D$ is an open continuous
surjection and $U\subseteq L(K)$ is open and locally connected, then $p(U)$ is
open and locally connected, so $L(D)$ is dense whenever $L(K)$ is. Thus
Theorem~\ref{t:open-map} cannot produce, from Theorems~\ref{t:arc-dense}
and~\ref{t:fat}, a continuum in $\DCP$ with $L(K)$ not dense. Such continua
are given below (Theorems~\ref{t:combfan} and~\ref{t:grid}).
\end{remark}

\subsection*{The Cantor comb, the Cantor fan and cones}\label{s:combfan}

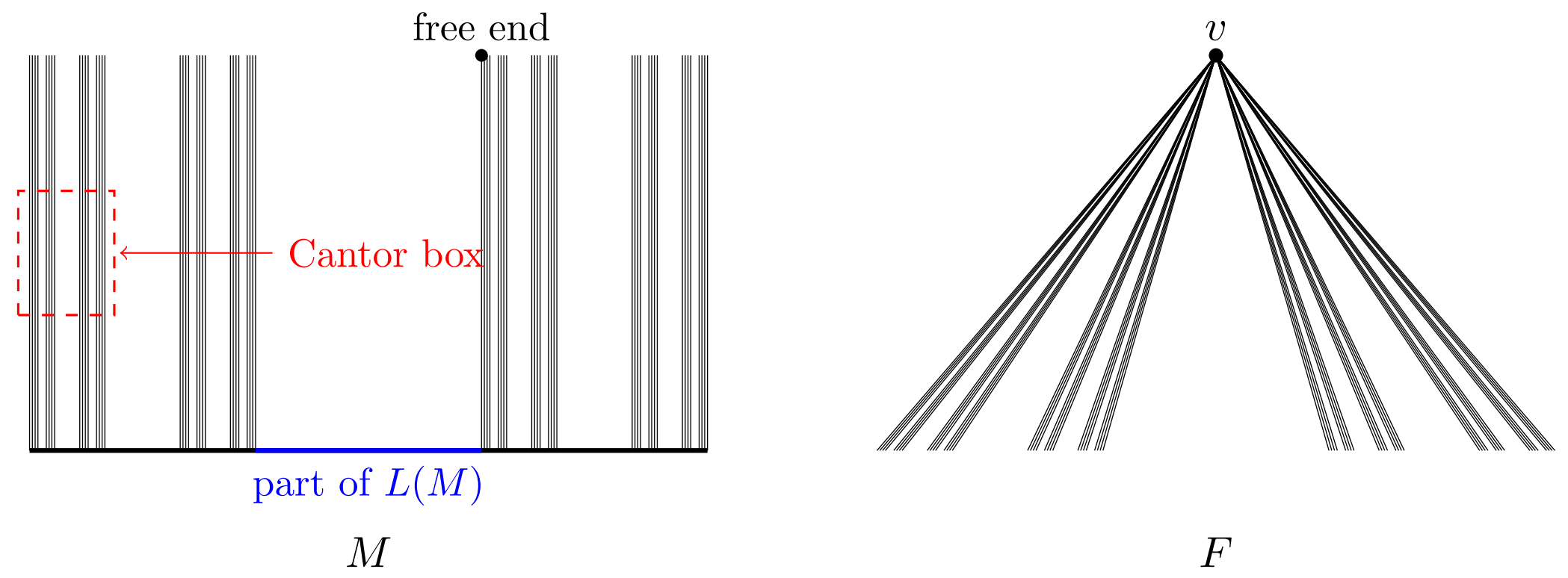
\begin{figure}[ht]
\centering
\begin{tikzpicture}[teeth/.style={line width=0.25pt},bar/.style={line width=1.2pt}]
\draw[teeth] (0.0000,0)--(0.0000,3.5);
\draw[teeth] (0.0247,0)--(0.0247,3.5);
\draw[teeth] (0.0494,0)--(0.0494,3.5);
\draw[teeth] (0.0741,0)--(0.0741,3.5);
\draw[teeth] (0.1481,0)--(0.1481,3.5);
\draw[teeth] (0.1728,0)--(0.1728,3.5);
\draw[teeth] (0.1975,0)--(0.1975,3.5);
\draw[teeth] (0.2222,0)--(0.2222,3.5);
\draw[teeth] (0.4444,0)--(0.4444,3.5);
\draw[teeth] (0.4691,0)--(0.4691,3.5);
\draw[teeth] (0.4938,0)--(0.4938,3.5);
\draw[teeth] (0.5185,0)--(0.5185,3.5);
\draw[teeth] (0.5926,0)--(0.5926,3.5);
\draw[teeth] (0.6173,0)--(0.6173,3.5);
\draw[teeth] (0.6420,0)--(0.6420,3.5);
\draw[teeth] (0.6667,0)--(0.6667,3.5);
\draw[teeth] (1.3333,0)--(1.3333,3.5);
\draw[teeth] (1.3580,0)--(1.3580,3.5);
\draw[teeth] (1.3827,0)--(1.3827,3.5);
\draw[teeth] (1.4074,0)--(1.4074,3.5);
\draw[teeth] (1.4815,0)--(1.4815,3.5);
\draw[teeth] (1.5062,0)--(1.5062,3.5);
\draw[teeth] (1.5309,0)--(1.5309,3.5);
\draw[teeth] (1.5556,0)--(1.5556,3.5);
\draw[teeth] (1.7778,0)--(1.7778,3.5);
\draw[teeth] (1.8025,0)--(1.8025,3.5);
\draw[teeth] (1.8272,0)--(1.8272,3.5);
\draw[teeth] (1.8519,0)--(1.8519,3.5);
\draw[teeth] (1.9259,0)--(1.9259,3.5);
\draw[teeth] (1.9506,0)--(1.9506,3.5);
\draw[teeth] (1.9753,0)--(1.9753,3.5);
\draw[teeth] (2.0000,0)--(2.0000,3.5);
\draw[teeth] (4.0000,0)--(4.0000,3.5);
\draw[teeth] (4.0247,0)--(4.0247,3.5);
\draw[teeth] (4.0494,0)--(4.0494,3.5);
\draw[teeth] (4.0741,0)--(4.0741,3.5);
\draw[teeth] (4.1481,0)--(4.1481,3.5);
\draw[teeth] (4.1728,0)--(4.1728,3.5);
\draw[teeth] (4.1975,0)--(4.1975,3.5);
\draw[teeth] (4.2222,0)--(4.2222,3.5);
\draw[teeth] (4.4444,0)--(4.4444,3.5);
\draw[teeth] (4.4691,0)--(4.4691,3.5);
\draw[teeth] (4.4938,0)--(4.4938,3.5);
\draw[teeth] (4.5185,0)--(4.5185,3.5);
\draw[teeth] (4.5926,0)--(4.5926,3.5);
\draw[teeth] (4.6173,0)--(4.6173,3.5);
\draw[teeth] (4.6420,0)--(4.6420,3.5);
\draw[teeth] (4.6667,0)--(4.6667,3.5);
\draw[teeth] (5.3333,0)--(5.3333,3.5);
\draw[teeth] (5.3580,0)--(5.3580,3.5);
\draw[teeth] (5.3827,0)--(5.3827,3.5);
\draw[teeth] (5.4074,0)--(5.4074,3.5);
\draw[teeth] (5.4815,0)--(5.4815,3.5);
\draw[teeth] (5.5062,0)--(5.5062,3.5);
\draw[teeth] (5.5309,0)--(5.5309,3.5);
\draw[teeth] (5.5556,0)--(5.5556,3.5);
\draw[teeth] (5.7778,0)--(5.7778,3.5);
\draw[teeth] (5.8025,0)--(5.8025,3.5);
\draw[teeth] (5.8272,0)--(5.8272,3.5);
\draw[teeth] (5.8519,0)--(5.8519,3.5);
\draw[teeth] (5.9259,0)--(5.9259,3.5);
\draw[teeth] (5.9506,0)--(5.9506,3.5);
\draw[teeth] (5.9753,0)--(5.9753,3.5);
\draw[teeth] (6.0000,0)--(6.0000,3.5);
\draw[bar] (0,0)--(6,0);
\draw[bar,blue] (2,0)--(4,0);
\node[blue,below] at (3,0) {\small part of $L(M)$};
\draw[red,dashed,line width=0.6pt] (-0.1,1.2) rectangle (0.75,2.3);
\draw[red,->] (2.15,1.75)--(0.8,1.75);
\node[red,right] at (2.15,1.75) {\small Cantor box};
\fill (4,3.5) circle (1.6pt);
\node[above] at (4,3.5) {\small free end};
\node at (3,-0.9) {$M$};
\begin{scope}[xshift=7.5cm]
\draw[teeth] (3,3.5)--(0.0000,0);
\draw[teeth] (3,3.5)--(0.0247,0);
\draw[teeth] (3,3.5)--(0.0494,0);
\draw[teeth] (3,3.5)--(0.0741,0);
\draw[teeth] (3,3.5)--(0.1481,0);
\draw[teeth] (3,3.5)--(0.1728,0);
\draw[teeth] (3,3.5)--(0.1975,0);
\draw[teeth] (3,3.5)--(0.2222,0);
\draw[teeth] (3,3.5)--(0.4444,0);
\draw[teeth] (3,3.5)--(0.4691,0);
\draw[teeth] (3,3.5)--(0.4938,0);
\draw[teeth] (3,3.5)--(0.5185,0);
\draw[teeth] (3,3.5)--(0.5926,0);
\draw[teeth] (3,3.5)--(0.6173,0);
\draw[teeth] (3,3.5)--(0.6420,0);
\draw[teeth] (3,3.5)--(0.6667,0);
\draw[teeth] (3,3.5)--(1.3333,0);
\draw[teeth] (3,3.5)--(1.3580,0);
\draw[teeth] (3,3.5)--(1.3827,0);
\draw[teeth] (3,3.5)--(1.4074,0);
\draw[teeth] (3,3.5)--(1.4815,0);
\draw[teeth] (3,3.5)--(1.5062,0);
\draw[teeth] (3,3.5)--(1.5309,0);
\draw[teeth] (3,3.5)--(1.5556,0);
\draw[teeth] (3,3.5)--(1.7778,0);
\draw[teeth] (3,3.5)--(1.8025,0);
\draw[teeth] (3,3.5)--(1.8272,0);
\draw[teeth] (3,3.5)--(1.8519,0);
\draw[teeth] (3,3.5)--(1.9259,0);
\draw[teeth] (3,3.5)--(1.9506,0);
\draw[teeth] (3,3.5)--(1.9753,0);
\draw[teeth] (3,3.5)--(2.0000,0);
\draw[teeth] (3,3.5)--(4.0000,0);
\draw[teeth] (3,3.5)--(4.0247,0);
\draw[teeth] (3,3.5)--(4.0494,0);
\draw[teeth] (3,3.5)--(4.0741,0);
\draw[teeth] (3,3.5)--(4.1481,0);
\draw[teeth] (3,3.5)--(4.1728,0);
\draw[teeth] (3,3.5)--(4.1975,0);
\draw[teeth] (3,3.5)--(4.2222,0);
\draw[teeth] (3,3.5)--(4.4444,0);
\draw[teeth] (3,3.5)--(4.4691,0);
\draw[teeth] (3,3.5)--(4.4938,0);
\draw[teeth] (3,3.5)--(4.5185,0);
\draw[teeth] (3,3.5)--(4.5926,0);
\draw[teeth] (3,3.5)--(4.6173,0);
\draw[teeth] (3,3.5)--(4.6420,0);
\draw[teeth] (3,3.5)--(4.6667,0);
\draw[teeth] (3,3.5)--(5.3333,0);
\draw[teeth] (3,3.5)--(5.3580,0);
\draw[teeth] (3,3.5)--(5.3827,0);
\draw[teeth] (3,3.5)--(5.4074,0);
\draw[teeth] (3,3.5)--(5.4815,0);
\draw[teeth] (3,3.5)--(5.5062,0);
\draw[teeth] (3,3.5)--(5.5309,0);
\draw[teeth] (3,3.5)--(5.5556,0);
\draw[teeth] (3,3.5)--(5.7778,0);
\draw[teeth] (3,3.5)--(5.8025,0);
\draw[teeth] (3,3.5)--(5.8272,0);
\draw[teeth] (3,3.5)--(5.8519,0);
\draw[teeth] (3,3.5)--(5.9259,0);
\draw[teeth] (3,3.5)--(5.9506,0);
\draw[teeth] (3,3.5)--(5.9753,0);
\draw[teeth] (3,3.5)--(6.0000,0);
\fill (3,3.5) circle (1.8pt);
\node[above] at (3,3.5) {$v$};
\node at (3,-0.9) {$F$};
\end{scope}
\end{tikzpicture}
\caption{The Cantor comb $M$ and the Cantor fan $F$, with $C$ replaced by the $64$ end points of the fifth step of its construction. In $M$ the base over the gap $(\tfrac13,\tfrac23)$ (blue) is a part of $L(M)$; a Cantor box is outlined in red. The roots of the teeth are on the base of $M$ and at the vertex $v$ of $F$.}\label{f:combfan}
\end{figure}

All positive results above for arbitrary $f$ assume that $L(K)$ is dense, and
by Remark~\ref{r:Ldense} open maps cannot remove this assumption. We now show
that two classical continua for which $L(K)$ is not dense are in $\DCP$. The
\emph{Cantor comb} is
\[
M=(C\times[0,1])\cup([0,1]\times\{0\}),
\]
the comb of~\cite[Problem~10]{Ki13}. The \emph{Cantor fan} $F$ is the cone
$C\times[0,1]/C\times\{0\}$ over the Cantor set; we use its planar model, the
union of the straight line segments (\emph{legs}) joining the \emph{vertex}
$v=(\tfrac12,1)$ to the points $(x,0)$, $x\in C$ (Figure~\ref{f:combfan}). The
points of the base $[0,1]\times\{0\}$ of $M$ over the gaps of $C$ have open
arcs as neighbourhoods. Every other point of $M$, and every point of $F$, has
a Cantor box (see Lemma~\ref{l:teeth}(ii) below) in each of its
neighbourhoods; since an open subset of a locally connected space is locally
connected and a Cantor box is not, such a point has no locally connected open
neighbourhood. Hence $L(M)=((0,1)\setminus C)\times\{0\}$, which is not dense
in $M$, and $L(F)=\emptyset$.

The proof for the fan works for all cones. For a nonempty compact metrizable
space $Y$ let $\operatorname{Cone}(Y)=Y\times[0,1]/Y\times\{0\}$, with vertex
$v$, the image of $Y\times\{0\}$; thus $F=\operatorname{Cone}(C)$. We treat
the comb and the cones at once. Let $K$ be $M$, and then put $Y=C$, or
$K=\operatorname{Cone}(Y)$, and define $h:Y\times[0,1]\to K$ and the
\emph{root set} $\mathcal R\subseteq K$ by
\[
h(x,\tau)=(x,\tau),\quad \mathcal R=[0,1]\times\{0\}\quad\text{for }K=M;
\]
$h$ the quotient map and $\mathcal R=\{v\}$ for $K=\operatorname{Cone}(Y)$.
For $F$ in its planar model, $h(x,\tau)=(1-\tau)v+\tau(x,0)$. The arcs
$T_x=h(\{x\}\times[0,1])$, $x\in Y$, are the \emph{teeth} of $K$ (the legs of
$F$); $h(x,0)\in\mathcal R$ is the \emph{root} of $T_x$, and $h(x,1)$ is its
free end point.

\begin{lemma}\label{l:teeth}
\begin{enumerate}[(i)]
\item $K$ is an arc connected continuum, $h$ is continuous,
$K=h(Y\times[0,1])\cup\mathcal R$, and $h$ maps $Y\times(0,1]$
homeomorphically onto the open set $K\setminus\mathcal R$.
\item For $A$ open in $Y$ and $J$ open in $(0,1]$ the set $h(A\times J)$ is
open in $K$; if $A\ne\emptyset$ and $J\subseteq(0,1)$ is an open interval, it
is a box (Remark~\ref{r:box}), and a Cantor box if $Y=C$ and $A$ is clopen.
Every point of $h(Y\times(0,1))$ has arbitrarily small box neighbourhoods
$h(A\times(\tau-\varepsilon,\tau+\varepsilon))$.
\item Every point of $\mathcal R$ has arbitrarily small open neighbourhoods
$U'$ that are \emph{closed downwards}: if $h(x,\tau)\in U'$, then
$h(x,\sigma)\in U'$ for $0\le\sigma\le\tau$.
\end{enumerate}
\end{lemma}

\begin{proof}
(i) For $M$ this is clear. For $\operatorname{Cone}(Y)$, the restriction of
the quotient map $h$ to the saturated open set $Y\times(0,1]$ is a bijective
quotient map onto the open set $K\setminus\{v\}$, hence a homeomorphism;
every point is joined to $v$ by a tooth. (ii) follows from (i). (iii) For
$u=(x_0,0)\in M$ take
$U'=M\cap\bigl((x_0-\varepsilon,x_0+\varepsilon)\times[0,\varepsilon)\bigr)$;
for $\operatorname{Cone}(Y)$ take $U'=h(Y\times[0,\delta))$, whose complement
$h(Y\times[\delta,1])$ is compact; these sets form a base at $v$ because $Y$
is compact.
\end{proof}

For the rest of this subsection $f\in C(K)$ and an $f$-admissible $g$ are
fixed, $\Sigma$ is the set of points of $K$ whose level is collapsed, and
$\Omega=K\setminus\cl\Sigma$. Thus $\Omega$ is open, and an open set disjoint
from $\Sigma$ is contained in $\Omega$.

\begin{lemma}\label{l:tooth-omega}
$f$ is constant on no nonempty open subset of $\Omega$, and $f$ is constant on
every connected subset of $\Omega$ that is contained in a tooth.
\end{lemma}

\begin{proof}
A nonempty open set on which $f$ is constant lies in $\Sigma$ by
Lemma~\ref{l:interior}, so it is not contained in $\Omega$. Suppose $f$ is not
constant on a connected set $J\subseteq\Omega$ contained in a tooth $T_x$;
then $J=h(\{x\}\times I_J)$ for an interval $I_J\subseteq[0,1]$. Choose
$\tau_1<\tau_2$ in $I_J$ with $f(h(x,\tau_1))\ne f(h(x,\tau_2))$, and a level
$t$ strictly between these values such that $f(h(x,\cdot))$ is not identically
$t$ on any nondegenerate subinterval of $[\tau_1,\tau_2]$ (all but countably
many $t$ qualify). As in the proof of Lemma~\ref{l:baire} there is
$\tau\in(\tau_1,\tau_2)$ such that $s=h(x,\tau)$ satisfies $f(s)=t$ and every
subinterval of $[\tau_1,\tau_2]$ containing $\tau$ in its interior contains
points where $f(h(x,\cdot))<t$ and points where it is $>t$. Since $0<\tau<1$ and
$s\in J\subseteq\Omega$, by Lemma~\ref{l:teeth}(ii) there is a box
$D=h(A\times(\tau-\varepsilon,\tau+\varepsilon))\subseteq\Omega$, and $s$ is a
crossing point of $D$. By the first part $f$ is constant on no nonempty open
subset of $D$, so by Proposition~\ref{p:box} (and Remark~\ref{r:box}) the
level $t$ is collapsed, i.e.\
$s\in\Sigma\cap\Omega=\emptyset$, a contradiction.
\end{proof}

\begin{lemma}\label{l:root}
$\mathcal R\subseteq\cl\Sigma$. Consequently $\Omega\subseteq h(Y\times(0,1])$.
\end{lemma}

\begin{proof}
\emph{The comb.} Let $V$ be an open interval of the base over a gap of $C$; it
is open in $M$, connected and locally connected. If $f$ is constant on a
nonempty open subset of $V$, that subset lies in $\Sigma$ by
Lemma~\ref{l:interior}. Otherwise $f(V)$ is a nondegenerate interval and, as
in the proof of Lemma~\ref{l:dense}, $V$ contains an openness point of $f$,
whose level is collapsed by Lemma~\ref{l:lc} (with $W=V$). So $\Sigma$ meets
every such $V$, and $\mathcal R=\cl L(M)\subseteq\cl\Sigma$.

\emph{The cones.} Suppose $v\in\Omega$. By Lemma~\ref{l:teeth}(iii) there is
$\delta>0$ such that the open set $S=h(Y\times[0,\delta))$ is contained in
$\Omega$. For each $x\in Y$ the set $h(\{x\}\times[0,\delta))$ is connected,
contained in $\Omega$ and in the tooth $T_x$, and contains $v$; by
Lemma~\ref{l:tooth-omega}, $f\equiv f(v)$ on it. Hence $f$ is constant on the
nonempty open set $S\subseteq\Omega$, contrary to Lemma~\ref{l:tooth-omega}.
So $v\in\cl\Sigma$.

Finally $K\setminus\mathcal R=h(Y\times(0,1])$.
\end{proof}

\begin{lemma}\label{l:tooth-approx}
Let $c\in f(K)$. Every $u\in Z_c^+$ admits a sequence $y_n$ as in
condition~(a) of Theorem~\ref{t:crit}. The same holds for $Z_c^-$.
\end{lemma}

\begin{proof}
The statement for $Z_c^-$ follows by applying the statement for $Z_{-c}^+$
to $-f$, which has the same fibres, the same $M_f$, the same admissible
functions and the same $\Sigma$.

If $u\in\cl(\Sigma\cap\{f>c\})$, take $y_n\in\Sigma\cap\{f>c\}$ with
$y_n\to u$. Otherwise there is an open neighbourhood $U$ of $u$ with
$U\cap\{f>c\}\subseteq\Omega$; then $U'\cap\{f>c\}$ is a nonempty open subset
of $\Omega$ for every open $U'\subseteq U$ containing $u$, because
$u\in\cl\{f>c\}$.

\emph{$u$ is not in $\mathcal R$.} Otherwise choose an open neighbourhood
$U'\subseteq U$ of $u$ that is closed downwards (Lemma~\ref{l:teeth}(iii)),
and put $W'=U'\cap\{f>c\}$. By Lemma~\ref{l:root}, $W'\cap\mathcal
R=\emptyset$; take $p=h(x,\tau)\in W'$, $\tau>0$. The set
$\{\sigma\in(0,1]:h(x,\sigma)\in W'\}$ is open in $(0,1]$; let $I$ be its
component containing $\tau$, and $\alpha=\inf I<\tau$. By
Lemma~\ref{l:tooth-omega}, $f\equiv f(p)>c$ on $h(\{x\}\times\cl I)$. The
point $e=h(x,\alpha)$ lies in $U'$, since $U'$ is closed downwards, and
$f(e)>c$, so $e\in W'$. If $\alpha=0$, then $e\in\mathcal R$, contrary to
$W'\cap\mathcal R=\emptyset$. If $\alpha>0$, then an interval about $\alpha$
lies in the open set $\{\sigma:h(x,\sigma)\in W'\}$, contrary to
$\alpha=\inf I$.

\emph{$u=h(x_u,\tau_u)$ with $0<\tau_u\le1$} (this includes the free end
points). Choose $0<\varepsilon<\tau_u$ and an open $A_u\ni x_u$ in $Y$ such
that $U_0=h(A_u\times J_u)\subseteq U$, where
$J_u=(\tau_u-\varepsilon,\tau_u+\varepsilon)\cap(0,1]$; $U_0$ is open by
Lemma~\ref{l:teeth}(ii). Put $W=U_0\cap\{f>c\}$, a nonempty open subset of
$\Omega$ with $u\in\cl W$. Let $x\in A_u$ and let $I$ be a component of the
open set $\{\sigma\in J_u:h(x,\sigma)\in W\}$. By Lemma~\ref{l:tooth-omega},
$f$ is constant, $>c$, on $h(\{x\}\times\cl I)$. A boundary point $\sigma$ of
$I$ in $J_u$ would therefore satisfy $h(x,\sigma)\in W$, and an interval about
$\sigma$ would lie in $\{\sigma'\in J_u:h(x,\sigma')\in W\}$, contrary to the
maximality of $I$; so $I=J_u$. Hence $W=h(A_W\times J_u)$, where
$A_W=\{x\in A_u:h(x,\tau_u)\in W\}$ is open in $Y$, and
$f(h(x,\sigma))=\varphi(x):=f(h(x,\tau_u))>c$ for $x\in A_W$, $\sigma\in J_u$.
Moreover $x_u\in\cl A_W$, and $\varphi$, defined by the same formula on $A_u$,
is continuous with $\varphi(x_u)=c$. For $x\in A_W$ put
\[
\beta(x)=\inf\{s\in[0,\tau_u]:\ f(h(x,\sigma))=\varphi(x)\text{ for all }
\sigma\in[s,\tau_u]\},\qquad e_x=h(x,\beta(x)).
\]
Then $\beta(x)\le\tau_u-\varepsilon$, and $f(h(x,\sigma))=\varphi(x)$ for
$\sigma\in[\beta(x),\tau_u]$; in particular $f(e_x)=\varphi(x)>c$. We have
$e_x\in\cl\Sigma$: if $\beta(x)=0$, then $e_x\in\mathcal R\subseteq\cl\Sigma$
by Lemma~\ref{l:root}; if $\beta(x)>0$ and $e_x\in\Omega$, then by
Lemma~\ref{l:tooth-omega} $f\equiv f(e_x)=\varphi(x)$ on a neighbourhood of
$e_x$ in $T_x$, contrary to the definition of $\beta(x)$. The function
$\beta$ is lower semicontinuous on $A_W$, because for $0\le s\le\tau_u$ the
set $\{x\in A_W:\beta(x)\le s\}=\{x\in A_W:f(h(x,\sigma))=\varphi(x)\text{ for
all }\sigma\in[s,\tau_u]\}$ is closed in $A_W$. Since $A_W$ is open in the
compact metrizable space $Y$, it is a Baire space, and the points of
continuity of $\beta$ are dense in $A_W$.

Let $x\in A_W$ be a point of continuity of $\beta$, and $0<\eta<\varepsilon$.
There is an open neighbourhood $A''\subseteq A_W$ of $x$ with
$\beta<\beta(x)+\eta$ on $A''$. The set
$O=h\bigl(A''\times((\beta(x)+\eta,\tau_u+\varepsilon)\cap(0,1])\bigr)$ is
open and contains $h(x,\tau_u)$, and $f(h(x',\sigma))=\varphi(x')$ on $O$: for
$\sigma\ge\tau_u$ because $h(x',\sigma)\in W$, and for $\sigma<\tau_u$ because
$\sigma>\beta(x')$. So $r(h(x',\sigma))=h(x',\beta(x)+\eta)$ defines a
continuous map $r:O\to K$ with $f\circ r=f$ on $O$, and by
Lemma~\ref{l:locret} $g(h(x,\tau_u))=g(h(x,\beta(x)+\eta))$. Letting
$\eta\to0$, $g(h(x,\tau_u))=g(e_x)$.

Choose points of continuity $x_n$ of $\beta$ with $x_n\to x_u$ and put
$y_n=e_{x_n}$. Then $y_n\in\cl\Sigma$, $f(y_n)=\varphi(x_n)>c$,
$f(y_n)\to\varphi(x_u)=c$ and $g(y_n)=g(h(x_n,\tau_u))\to g(u)$.
\end{proof}

\begin{theorem}\label{t:combfan}
The Cantor comb $M$ and the cone $\operatorname{Cone}(Y)$ over every nonempty
compact metrizable space $Y$ are in $\DCP$. In particular the Cantor fan $F$
is in $\DCP$, although $L(M)$ is not dense in $M$ and $L(F)=\emptyset$.
\end{theorem}

\begin{proof}
Let $K$ be $M$ or $\operatorname{Cone}(Y)$, $f\in C(K)$, $g$ $f$-admissible and $c\in f(K)$. If
$\Int Z_c\ne\emptyset$, the level $c$ is collapsed by Lemma~\ref{l:interior}.
Otherwise condition~(a) of Theorem~\ref{t:crit} holds by
Lemma~\ref{l:tooth-approx}, and (b) holds because $K$ is arc connected (see
the remark preceding Theorem~\ref{t:arc-dense}). Hence every level is
collapsed.
\end{proof}

\begin{corollary}\label{c:smooth}
Every smooth fan whose set of end points is closed is in $\DCP$.
\end{corollary}

\begin{proof}
Such a fan is either an arc, which is locally connected, or an open image of
$F$~\cite[Theorem~8]{CC89}; apply Theorem~\ref{t:open-map}.
\end{proof}

\begin{remark}\label{r:combfan}
(1) Theorem~\ref{t:combfan} answers Problem~10 of~\cite{Ki13} and the
question, raised there, whether the Cantor fan is in $\DCP$.

(2) The proof follows each plateau of $f$ on the arcs of a box towards the
root of the tooth; the free end points of the teeth are never needed. What is
needed is that the roots lie in $\cl\Sigma$ for every $f$: for the comb this
holds because $\mathcal R=\cl L(M)$, for the cones because of the argument at
the vertex in Lemma~\ref{l:root}.

(3) For cones, $Y$ need not be connected, nor in $\DCP$, nor contain any arc.
Besides the Cantor fan, Theorem~\ref{t:combfan} applies, e.g., to the cones over
solenoids, over the buckethandle continuum and over the pseudo-arc; for these
$L(Y)$ is not dense, so Theorem~\ref{t:arc-dense} does not apply. (Since
$\operatorname{Cone}(Y)\setminus\{v\}\cong Y\times(0,1]$, we have
$L(\operatorname{Cone}(Y))\setminus\{v\}=h(L(Y)\times(0,1])$, so
$L(\operatorname{Cone}(Y))$ is dense if and only if $L(Y)$ is. For the harmonic
fan $\operatorname{Cone}(\{0\}\cup\{1/n:n\in\N\})$, $L(Y)$ is dense and
Theorem~\ref{t:arc-dense} applies as well; compare Example~\ref{e:cone}.) By
Corollary~\ref{c:product}, products behave
differently: for the $\sin(1/x)$ continuum $S$ of Example~\ref{e:cone},
$\operatorname{Cone}(S)\in\DCP$ but $S\times[0,1]\notin\DCP$. By
Theorem~\ref{t:open-map}, open images of cones are in $\DCP$.

(4) The same proof, with $\mathcal R=[0,1]\times\{0,1\}$ and, at the points of
the top bar, neighbourhoods closed upwards, shows that the \emph{Cantor
ladder} $(C\times[0,1])\cup([0,1]\times\{0,1\})$ is in $\DCP$; similarly, with
two vertices, the suspension of any nonempty compact metrizable space is in
$\DCP$.
\end{remark}

\subsection*{Rooted teeth}\label{s:teeth}

The proofs of Lemmas~\ref{l:tooth-omega}--\ref{l:tooth-approx} use only the
following structure, which also underlies Theorems~\ref{t:arc-dense}
and~\ref{t:fat}.

Recall from Definition~\ref{d:tooth} the notions of a tooth chart
$h:U\times[0,1]\to K$ with $h(U\times J_h)=O_h$ open, of open teeth
$h(\{y\}\times J_h)$ and their roots $h(y,0)$, and of star points.

\begin{theorem}\label{t:teeth}
Let $K$ be a metrizable continuum with a family of tooth charts, let $B$ be
the union of the sets $O_h$, and let $\mathcal R\subseteq K$ be closed. Assume:
\begin{enumerate}[(1)]
\item $K$ is arc connected, or every arc component $P$ of $K$ is dense or
satisfies $P\subseteq\cl\Int P$;
\item $K=\cl L(K)\cup\mathcal R\cup B$, and $h(U\times\{0\})\subseteq\mathcal
R\cup\cl L(K)$ for every chart $h$;
\item every point of $(\mathcal R\cup\cl L(K))\setminus(L(K)\cup B)$ has
arbitrarily small open neighbourhoods $U'$ such that every component of
$U'\cap T$, $T$ an open tooth, has an end point in $U'\setminus B$ (an end
point of $h(\{y\}\times I)$, $I$ an interval, is a point $h(y,\alpha)$ with
$\alpha$ an end point of $I$);
\item every point of $\mathcal R$ lies in $\cl L(K)$ or is a star point.
\end{enumerate}
Then $K\in\DCP$.
\end{theorem}

\begin{proof}
Fix $f\in C(K)$ and an $f$-admissible $g$, and let $\Sigma$ and
$\Omega=K\setminus\cl\Sigma$ be as before.

(i) \emph{$L(K)\subseteq\cl\Sigma$}: every nonempty open subset of $L(K)$
contains a nonempty connected open locally connected set, which contains a
point of $D_f\subseteq\Sigma$ by the proof of Lemma~\ref{l:dense}.

(ii) As in Lemma~\ref{l:tooth-omega}, $f$ is constant on no nonempty open
subset of $\Omega$ and on every connected subset of $\Omega$ contained in an
open tooth: the crossing point $s=h(y,\tau)$ found there has $0<\tau<1$, so it
has a box neighbourhood $h(U''\times(\tau-\delta,\tau+\delta))\subseteq\Omega$,
and Proposition~\ref{p:box} with Remark~\ref{r:box} applies.

(iii) \emph{$\mathcal R\subseteq\cl\Sigma$, hence $\Omega\subseteq B$}: points
of $\mathcal R\cap\cl L(K)$ are in $\cl\Sigma$ by (i); if a star point $v$ were
in $\Omega$, then, as for the vertex in Lemma~\ref{l:root}, $f$ would be
constant, equal to $f(v)$, on a star neighbourhood $S\subseteq\Omega$, contrary
to (ii). By (2), $\Omega\subseteq B$.

(iv) \emph{Condition (a) of Theorem~\ref{t:crit}.} Let $u\in Z_c^+$ ($Z_c^-$
is handled by passing to $-f$). If $u\notin\cl(\Sigma\cap\{f>c\})$, fix an
open $N\ni u$ with $N\cap\{f>c\}\subseteq\Omega\subseteq B$. Since $u\in
\cl(N\cap\{f>c\})\subseteq\cl\Omega$ and $L(K)$ is open and disjoint from
$\Omega$, $u\notin L(K)$. If $u\in B$, then $u=h(x_u,\tau_u)$ with $\tau_u\in
J_h$, and the second case of the proof of Lemma~\ref{l:tooth-approx} applies in
the chart $h$, with $J_u=(\tau_u-\varepsilon,\tau_u+\varepsilon)\cap J_h$. (There
the openness of the sets $h(A\times J)$ now follows from $O_h$ being open and
$h$ being a homeomorphism on $U\times J_h$; Baire's theorem is applied to an
open subset of $U$; and the continuity of $h$ on $U\times[0,1]$ makes the
relevant sublevel sets of $\beta$ closed and puts each $e_x$ in the closure of
its plateau.) The end points $e_x$ of the plateaus lie in $\cl\Sigma$, either
as roots, which lie in $\mathcal R\cup \cl L(K)\subseteq\cl\Sigma$, or by (ii).
Otherwise $u\in(\mathcal R\cup\cl L(K))\setminus(L(K)\cup B)$ by (2). Take
$U'\subseteq N$ as in (3), put
$W'=U'\cap\{f>c\}$, a nonempty open subset of $\Omega$, and let $p\in W'$ lie
on the open tooth $T$. Let $J$ be the component of $U'\cap T$ and $I\subseteq
J$ the component of $W'\cap T$ containing $p$. By (ii), $f\equiv f(p)>c$ on
$\cl I$, so a boundary point of $I$ in $J$ would lie in $W'$; hence $I=J$. The
end point $e\in U'\setminus B$ of $J$ given by (3) satisfies $f(e)>c$, so
$e\in W'\subseteq B$, a contradiction.

Finally, levels with interior are collapsed by Lemma~\ref{l:interior}, and
condition (b) holds by (1) (the remark preceding Theorem~\ref{t:arc-dense},
or Lemma~\ref{l:dense-arc}). Theorem~\ref{t:crit} applies to every level.
\end{proof}

\begin{remark}\label{r:teeth}
(1) Without charts ($B=\emptyset$) and with $\mathcal R=K\setminus L(K)$,
Theorem~\ref{t:teeth} is Theorems~\ref{t:arc-dense} and~\ref{t:fat}. The comb
and the cones are covered with the charts of Lemma~\ref{l:teeth}, with
$\mathcal R=C\times\{0\}\subseteq\cl L(M)$, resp.\ $\mathcal R=\{v\}$, a star
point; the Cantor ladder with $\mathcal R=C\times\{0,1\}\subseteq\cl L$. In a
cone, points of $\cl L\setminus L$ may lie in $B$ (e.g.\ the leg over $0$ in
the harmonic fan); this is why (3) is required only outside $B$. The
suspension of a nonempty compact metrizable space is covered with $\mathcal
R$ the two vertices, both star points.

(2) The Cantor grid below satisfies (1)--(3), with the nodes as $\mathcal R$
and the charts $(x,\tau)\mapsto(x,a+\tau(b-a))$ over the gaps $(a,b)$ (and
their horizontal analogues), but the nodes are neither in $\cl L(G)=\emptyset$
nor star points. In fact no family of tooth charts makes $G$ satisfy
Theorem~\ref{t:teeth}. Since $G$ is one-dimensional, each $O_h\cong U\times
J_h$ has $\dim U=0$, so $O_h$ contains no simple triod; as every neighbourhood
of a node contains one, no node lies in $B$. If $v\in G$ is not a node, every
sufficiently small neighbourhood $S$ of $v$ meets other segments parallel to
the one containing $v$, and their points cannot be joined to $v$ by an arc in
$S$; and every neighbourhood of a node contains other nodes, which lie in no
open tooth. Hence $G$ has no star points. Thus (4) forces $\mathcal
R\subseteq\cl L(G)=\emptyset$, and then (2) fails at the nodes. (Heuristically,
along a line through a node $f$ may vary like the Cantor function.) The main
step of the proof for the grid, Lemma~\ref{l:nodes}, shows
that nevertheless $\Omega$ contains no node; with this in place of (4) the
proof above becomes the proof of Theorem~\ref{t:grid}.
\end{remark}

\subsection*{The Cantor grid}\label{s:grid}

In the comb and the cones every arc of a box lies in a tooth whose root
is in $\cl\Sigma$ for every $f$ (Lemma~\ref{l:root}). We now give a continuum
in $\DCP$ in which, as in the fan, no point has an open locally connected
neighbourhood, but in which both ends of every such arc are points of a
Cantor set of \emph{nodes}, about which nothing of this kind is known a
priori.

Recall that $C$ is the ternary Cantor set. Let $\mathcal E$ be the set
of end points of the bounded components (\emph{gaps}) of $[0,1]\setminus C$,
and $\vartheta$ the Cantor function; thus $\vartheta(\mathcal E)$ is the set
$\Delta$ of dyadic rationals in $(0,1)$. (We write $\vartheta$ rather than
the usual $F$, which denotes the Cantor fan.) The
\emph{Cantor grid} is
\[
G=(C\times[0,1])\cup([0,1]\times C).
\]
The points of $C\times C$ are the \emph{nodes} of $G$. For a gap $(a,b)$ and
$x,y\in C$ the arcs $\{x\}\times[a,b]$ and $[a,b]\times\{y\}$ are the
(vertical and horizontal) \emph{segments} of $G$. Every point of $G$ that is not
a node lies in the relative interior of exactly one segment, and the sets
$(C\cap I)\times(a,b)$ and $(a,b)\times(C\cap I)$, $I$ an open interval with
end points not in $C$, are
Cantor boxes in $G$ (with vertical, resp.\ horizontal, arcs). Put
$\kappa=\vartheta \times \vartheta :G\to[0,1]^2$.

\begin{lemma}\label{l:grid}
\begin{enumerate}[(i)]
\item $G$ is an arc connected planar continuum and $L(G)=\emptyset$.
\item Let $z=(x,y)$ with $x,y\in C\setminus(\mathcal E\cup\{0,1\})$. Then $\kappa$ maps
every neighbourhood of $z$ in $G$ onto a neighbourhood of $\kappa(z)$.
\item If $\sigma,\tau\in(0,1)\setminus\Delta$, then $\kappa^{-1}(\sigma,\tau)$
is a single node, and its coordinates lie in $C\setminus \mathcal E$.
\end{enumerate}
\end{lemma}

\begin{proof}
(i) Every vertical line $\{x\}\times[0,1]$, $x\in C$, meets every horizontal
line $[0,1]\times\{y\}$, $y\in C$, and every point of $G$ lies on one of these
arcs. Every nonempty open subset of $G$ contains a Cantor box, whose
components are arcs that are not open in it; so no point has a locally
connected open neighbourhood. On the other hand, $G$ is locally connected at
every node $(x,y)$: for $\varepsilon>0$ the set
$G\cap\bigl((x-\varepsilon,x+\varepsilon)\times(y-\varepsilon,y+\varepsilon)\bigr)$
is arc connected, since each of its points is joined to $(x,y)$ within it by
a vertical and a horizontal straight line interval contained in $G$.

(ii) There are $x'<x<x''$ in $C$, arbitrarily close to $x$, with
$\vartheta (x')<\vartheta (x)<\vartheta (x'')$, because $x\in C$ is not an end point of a gap; similarly
$y'<y<y''$. If $\vartheta (x')<\sigma<\vartheta (x'')$ and $\vartheta (y')<\tau<\vartheta (y'')$, choose
$\xi\in C\cap[x',x'']$ with $\vartheta (\xi)=\sigma$ (possible because $\vartheta $ is constant on
the gaps) and $\eta\in[y',y'']$ with $\vartheta (\eta)=\tau$; then $(\xi,\eta)\in
C\times[0,1]\subseteq G$ and $\kappa(\xi,\eta)=(\sigma,\tau)$.

(iii) $\vartheta ^{-1}(\sigma)$ is either a point of $C\setminus \mathcal E$ or the closure of a
gap, and the latter happens exactly when $\sigma\in\Delta$.
\end{proof}

For the rest of this subsection $f\in C(G)$ and an $f$-admissible $g$ are fixed,
$\Sigma$ is the set of points of $G$ whose level is collapsed, and
$\Omega=G\setminus\cl\Sigma$. Thus $\Omega$ is open, and an open set disjoint
from $\Sigma$ is contained in $\Omega$.

\begin{lemma}\label{l:omega}
$f$ is constant on no nonempty open subset of $\Omega$, and $f$ is constant on
every connected subset of $\Omega$ that is contained in a segment.
\end{lemma}

\begin{proof}
The proof of Lemma~\ref{l:tooth-omega} applies, with the tooth replaced by a
segment, say $\{x\}\times[a,b]$, and with $\tau_1,\tau_2\in(a,b)$. The only
change is in the choice of the Cantor box: the crossing point $s$ lies in the
relative interior of the segment, so it is not a node and has a Cantor box
neighbourhood $B\subseteq\Omega$ with vertical arcs.
\end{proof}

\begin{lemma}\label{l:nodes}
$\Omega$ contains no node.
\end{lemma}

\begin{proof}
Suppose it does. Since the nodes with both coordinates in
$C\setminus(\mathcal E\cup\{0,1\})$ are dense in $C\times C$, there are $x_1<x_2$ and
$y_1<y_2$ in $C\setminus \mathcal E$ such that $S=G\cap([x_1,x_2]\times[y_1,y_2])$ is
contained in $\Omega$ and contains such a node in its interior. Since
$x_i,y_i\notin \mathcal E$, a gap meeting $(x_1,x_2)$ is contained in it, and similarly
for $(y_1,y_2)$. Put $Q=[\vartheta (x_1),\vartheta (x_2)]\times[\vartheta (y_1),\vartheta (y_2)]$; then
$\kappa(S)=Q$.

\emph{Step 1: $f=\psi\circ\kappa$ on $S$ for some $\psi\in C(Q)$.} A fibre of
$\kappa|_S$ is a node, a segment, or the union of the four segments bounding a
rectangle $[a,b]\times[a',b']$ with $(a,b)$, $(a',b')$ gaps. By
Lemma~\ref{l:omega} and continuity $f$ is constant on every segment contained
in $S$, hence on every fibre of $\kappa|_S$. Since $\kappa|_S$ is a closed
surjection onto $Q$, $\psi$ is well defined and continuous.

\emph{Step 2: $\psi$ is constant on no nonempty open subset of $Q$, and on no
nondegenerate segment $\{\sigma\}\times J\subseteq\Int Q$ or
$J\times\{\tau\}\subseteq\Int Q$ with $\sigma\in\Delta$, resp.\ $\tau\in\Delta$.}
If $\psi$ is constant on a nonempty open set $V\subseteq Q$, we may assume
$V\subseteq\Int Q$; then $\kappa^{-1}(V)$ is open in $G$, nonempty (since
$\vartheta (C)=[0,1]$) and contained in $S$ (a point of $G$ whose image lies in $\Int Q$
has $x_1<x<x_2$ and $y_1<y<y_2$), and $f$ is constant on it, contrary to
Lemma~\ref{l:omega}. If $\psi$ is constant on $\{\sigma\}\times J$ with
$\sigma=\vartheta (a)=\vartheta (b)$, $(a,b)$ a gap in $(x_1,x_2)$, then $f$ is constant on
$(a,b)\times(C\cap \vartheta ^{-1}(\Int J))$, a nonempty open subset of $G$ contained in
$S\subseteq\Omega$, contrary to Lemma~\ref{l:omega}. The horizontal case is
symmetric.

\emph{Step 3.} Let $O_\psi$ be the set of points $q\in\Int Q$ such that $\psi$
maps every neighbourhood of $q$ onto a neighbourhood of $\psi(q)$, and let $R$
be the set of nodes $z\in S$ with coordinates in $C\setminus \mathcal E$ and
$\kappa(z)\in O_\psi$. We claim:

$(\ast)$ \emph{if $V'\subseteq\Int Q$ is an open disc and $t$ is an interior
point of $\psi(V')$, then there is $z\in R$ with $\kappa(z)\in V'$ and
$f(z)=t$.}

Let $A=V'\cap\{\psi<t\}$ and $B=V'\cap\{\psi>t\}$; both are nonempty and open.
Put $M=V'\cap\cl A\cap\cl B$. A point of $V'\setminus(\cl A\cup\cl B)$ would
have a neighbourhood on which $\psi\equiv t$, which is excluded by Step~2.
Hence $V'\setminus M$ is the union of the disjoint open sets
$V'\setminus\cl B\supseteq A$ and $V'\setminus\cl A\supseteq B$, so the
relatively closed set $M$ separates $V'$ between a point of $A$ and a point of
$B$. Since $V'$ is homeomorphic to $\R^2$, it is connected, locally connected
and contractible with respect to the circle~\cite[\S57, I, Theorem~9(i)]{Kur68};
hence $M$ contains a set $K_0$, closed in $V'$ and connected, which separates
these points~\cite[\S57, III, Theorem~1]{Kur68}. In particular $K_0$ is not a
point. Let $p\in K_0$ and let $D\subseteq V'$ be a closed disc about $p$ with
$K_0\not\subseteq D$. Then $X=\cl K_0$ (closure in $\R^2$) is a continuum and
$X\cap D=K_0\cap D$ is a nonempty proper closed subset of $X$. By the boundary
bumping theorem~\cite[Theorem~5.7]{Na92} the component $K_1$ of $K_0\cap D$
containing $p$ meets $\cl(X\setminus D)\subseteq X\setminus\Int D$; so $K_1$
is a nondegenerate continuum contained in $K_0\cap D$.

Let $\Lambda$ be the union of the lines $\{\sigma\}\times\R$ and
$\R\times\{\tau\}$, $\sigma,\tau\in\Delta$. If $K_1\subseteq\Lambda$, then by
Baire's theorem one of the countably many closed sets $K_1\cap \ell$, $\ell$ such a
line, has nonempty interior $U$ in $K_1$. If $\cl U=K_1$, then $K_1\subseteq \ell$;
otherwise, by~\cite[Theorem~5.7]{Na92}, the component of $\cl U$ containing a
point of $U$ meets $\cl(K_1\setminus\cl U)$, so it is a nondegenerate
continuum contained in $\cl U\subseteq \ell$. Either way $M\subseteq\{\psi=t\}$ contains a nondegenerate segment of a line in
$\Lambda$, contrary to Step~2. Hence there is $q'\in K_1$ with both
coordinates outside $\Delta$. Every disc about $q'$ is connected and meets $A$
and $B$, so $q'\in O_\psi$. By Lemma~\ref{l:grid}(iii) $\kappa^{-1}(q')$ is a
single node $z$ with coordinates in $C\setminus \mathcal E$; it lies in $S$ because
$q'\in\Int Q$. Thus $z\in R$, $\kappa(z)=q'\in V'$ and $f(z)=\psi(q')=t$,
which proves $(\ast)$.

By Step~2 and $(\ast)$, $R\ne\emptyset$. Moreover $f|_R$ is open: let $s\in R$
and let $V$ be a neighbourhood of $s$ in $G$. The coordinates of $s$ are
different from $0$ and $1$ because $\kappa(s)\in\Int Q$, so by
Lemma~\ref{l:grid}(ii) there
is an open disc $V'\subseteq\kappa(V)\cap\Int Q$ about $\kappa(s)$, and since
$\kappa(s)\in O_\psi$, $\psi(V')$ contains an interval $(f(s)-\delta,
f(s)+\delta)$. For $t$ in this interval $(\ast)$ gives $z\in R$ with $f(z)=t$
and $\kappa(z)\in V'\subseteq\kappa(V)$; since $z$ is the only point of $G$
mapped to $\kappa(z)$, $z\in V$. By Lemma~\ref{l:open} every level in $f(R)$
is collapsed, i.e.\ $R\subseteq\Sigma$. But $R\subseteq S\subseteq\Omega$, a
contradiction.
\end{proof}

\begin{lemma}\label{l:approx}
Let $c\in f(G)$. Every $u\in Z_c^+$ admits a sequence $y_n$ as in
condition (a) of Theorem~\ref{t:crit}. The same holds for $Z_c^-$.
\end{lemma}

\begin{proof}
If $u\in\cl(\Sigma\cap\{f>c\})$, take $y_n\in\Sigma\cap\{f>c\}$ with $y_n\to u$.
Otherwise there is an open neighbourhood $U$ of $u$ with
$U\cap\{f>c\}\subseteq\Omega$.

\emph{$u$ is not a node.} Otherwise let $U'\subseteq U$ be the trace on $G$ of
an open square $(x_0-\varepsilon,x_0+\varepsilon)\times
(y_0-\varepsilon,y_0+\varepsilon)$ about $u=(x_0,y_0)$, and
$W=U'\cap\{f>c\}$, a nonempty open subset of $\Omega$. By
Lemma~\ref{l:nodes} $W$ contains no node, so a point $p\in W$ lies, say, in the
interior of a vertical segment $\{x\}\times[a,b]$. Let $I$ be the component of
$W\cap(\{x\}\times(a,b))$ containing $p$. By Lemma~\ref{l:omega} $f\equiv f(p)>c$
on $\cl I$. An end point of $I$ lying in $U'$ would therefore belong to $W$,
which is impossible if it lies in $\{x\}\times(a,b)$ (maximality of $I$) or is
a node (Lemma~\ref{l:nodes}). Hence $I=\{x\}\times(y_0-\varepsilon,
y_0+\varepsilon)$, so $y_0\in(a,b)$, contradicting $y_0\in C$. For a horizontal
segment we obtain $x_0\in(a,b)$ in the same way.

\emph{$u$ lies in a segment.} Say $u=(x_u,\tau_u)$ with $x_u\in C$ and
$\tau_u$ in a gap $(a,b)$. The proof of the case $u=h(x_u,\tau_u)$ of
Lemma~\ref{l:tooth-approx} applies, with Lemma~\ref{l:omega} in place of
Lemma~\ref{l:tooth-omega}, the tooth $T_x$ replaced by the segment
$\{x\}\times[a,b]$ (and $h(x,\sigma)$ by $(x,\sigma)$), and the following
changes: the Cantor box $A_u\times(\tau_u-\varepsilon,\tau_u+\varepsilon)$ is
chosen with $[\tau_u-\varepsilon,\tau_u+\varepsilon]\subseteq(a,b)$; in the
definition of $\beta$ the interval $[0,\tau_u]$ is replaced by $[a,\tau_u]$;
and if $\beta(x)=a$, then $e_x$ is a node, so $e_x\notin\Omega$ by
Lemma~\ref{l:nodes}. Points of horizontal segments are treated in the same
way, and $Z_c^-$ is symmetric.
\end{proof}

\begin{theorem}\label{t:grid}
The Cantor grid $G$ is in $\DCP$, although $L(G)=\emptyset$.
\end{theorem}

\begin{proof}
Let $f\in C(G)$, let $g$ be $f$-admissible and $c\in f(G)$. If
$\Int Z_c\ne\emptyset$, the level $c$ is collapsed by Lemma~\ref{l:interior}.
Otherwise condition (a) of Theorem~\ref{t:crit} holds by
Lemma~\ref{l:approx}, and (b) holds because $G$ is arc connected (see the
remark preceding Theorem~\ref{t:arc-dense}). Hence every level is
collapsed.
\end{proof}

\begin{remark}\label{r:grid}
(1) Every arc of a box of $G$ ends at nodes, where $G$ is locally connected
(although no neighbourhood of a node is locally connected, so this does not
contradict $L(G)=\emptyset$). Unlike the roots of the comb and the fan
(Lemma~\ref{l:root}), the nodes are not known a priori to lie outside
$\Omega$; Lemma~\ref{l:nodes} shows this, using that near a node a function on
which our lemmas give no information factors through the monotone map
$\kappa$ onto a square (Step~1). Lemma~\ref{l:locret} then carries the
information from the nodes along the arcs.

(2) In Lemmas~\ref{l:tooth-approx} and~\ref{l:approx} the points $y_n$
attached to a point $u$ of a box need not lie on collapsed levels themselves; this is why condition (a) of
Theorem~\ref{t:crit} is formulated with $y_n\in\cl\Sigma$.
\end{remark}

\subsection*{Necessary conditions and counterexamples}

\begin{proposition}[{cf.~\cite[Theorem~1.15]{Ki10}}]\label{p:conn}
Let $K$ be a compact metrizable space and let $u,v\in K$ be non-isolated points
lying in different components of $K$. Then $\D(K)\ne C(K)$. In particular a compact metrizable space without isolated
points is in $\DCP$ only if it is connected.
\end{proposition}

\begin{proof}
Components of a compact space coincide with quasi-components, so there is a
clopen $W$ with $u\in W$, $v\notin W$. Let $f=d(\cdot,u)$ on $W$ and
$f=-d(\cdot,v)$ on $K\setminus W$. Every fibre other than $Z_0=\{u,v\}$ lies in
$W$ or in $K\setminus W$. Since $u$ is the limit of points $k_n\in W\setminus\{u\}$,
and $M_f(k_n)\subseteq W$, lower semicontinuity gives $M_f(u)\subseteq W$;
similarly $M_f(v)\subseteq K\setminus W$. So $g=\chi_W$ is constant on every
$M_f(k)$ but not on $Z_0$, so $f\notin\D(K)$.
\end{proof}

\begin{remark}
In~\cite[Theorem~1.15]{Ki10} the corresponding argument uses the characteristic
function of a component, which need not be continuous; the clopen set above
repairs this and also covers components that are singletons.
\end{remark}

\begin{example}[{\cite[Example~1.16]{Ki10}}]\label{e:sin}
Let $S$ be the closure of $\{(x,\sin(1/x)):0<x\le1\}$ and $f(x,y)=x$. All fibres
except $E=\{0\}\times[-1,1]$ are singletons. If $k\in E$ and $k_n\to k$ are points
of the curve, then $M_f(k_n)=\{k_n\}$, so $M_f(k)=\{k\}$. Hence every
$g\in C(S)$ is $f$-admissible, and $f\notin\D(S)$.
\end{example}

The same argument applies to every compactification of a ray, and of two
disjoint rays each of which is open in it.

\begin{proposition}\label{p:ray}
Let $K$ be a continuum containing a dense subset $U$ which is the
disjoint union of one or two open subsets $R_1,R_2$ of $K$, each homeomorphic
to $[0,\infty)$, and suppose that $A=K\setminus U$ is not a single point. Then
$K\notin\DCP$.
\end{proposition}

\begin{proof}
Let $r_i:[0,\infty)\to R_i$ be homeomorphisms and define $f=0$ on $A$,
$f(r_1(t))=1/(1+t)$ and $f(r_2(t))=-1/(1+t)$. Since $R_1$ and $R_2$ are open,
$f$ is continuous on $U$; and a point $a\in A$ has, for every $T>0$, a neighbourhood
disjoint from the compact set $r_1([0,T])\cup r_2([0,T])$, on which
$|f|<1/(1+T)$. So $f$ is continuous, and its fibres are the points of $U$ and
the set $A$. Let $\Phi\in\mathfrak M_f$; then $\Phi(k)=\{k\}$ for $k\in U$.
If $a\in A$ and $b\in\Phi(a)$, $b\ne a$, take disjoint open sets $V\ni b$ and
$V'\ni a$; by lower semicontinuity some neighbourhood $N\subseteq V'$ of $a$
satisfies $\Phi(k)\cap V\ne\emptyset$ for $k\in N$, which fails for
$k\in N\cap U$. Hence $M_f(k)=\{k\}$ for all $k$, every $g\in C(K)$ is
$f$-admissible, and a $g$ that is not constant on $A$ shows $f\notin\D(K)$.
\end{proof}

\begin{remark}\label{r:ray}
Proposition~\ref{p:ray} applies to the $\sin(1/x)$ continuum, to a ray
spiralling onto a circle, to a simple triod with a ray spiralling around
it~\cite[Fig.~10.13]{Il25}, and to Nadler's spiral continua~\cite[Fig.~10.12]{Il25}.
It does not apply to a line both of whose ends accumulate on the same set $A$
(no continuous $f$ is then injective on the line), e.g.\ two $\sin(1/x)$
curves approaching one segment from both sides and joined to each other by an
arc; we do not know whether such continua are in $\DCP$.
\end{remark}

\begin{example}\label{e:sin-arc}
Let $K=S\cup A$, where $S$ is as in Example~\ref{e:sin} and $A=[-1,0]\times\{-1\}$,
and let $f(x,y)=x$. The set $L(K)$, which contains the curve and $A$ without
the point $(0,-1)$, is dense, and every point of $K\setminus L(K)=E$ can be joined
by an arc (along $E$ and $A$) to a point of $L(K)$. Nevertheless the fibres of
$f$ other than $E$ are singletons, and as in Example~\ref{e:sin} $M_f(k)=\{k\}$ for
all $k$, so $K\notin\DCP$. Thus in Theorem~\ref{t:arc-dense} arc connectedness
cannot be replaced by the condition that every point of $K$ can be joined by an
arc to $L(K)$ (condition (6) of~\cite[Proposition~7]{Ki13} without the other
conditions).
\end{example}

In Proposition~\ref{p:square} the size of the square matters.

\begin{example}\label{e:halfsquare}
With $C_2$ and $E$ as in Proposition~\ref{p:square}, let
$Q'=[-1,0]\times[0,1]$, let $S$ be the closure of $C_2$ (so
$S=C_2\cup E$), and $K=Q'\cup S$. The arc components are $C_2$ and $Q'\cup E$,
and $Q'\cup E\not\subseteq\cl\Int(Q'\cup E)=Q'$, while $L(K)\supseteq\Int
Q'\cup C_2$ is dense. Let $f(x,y)=x$. For $0<x\le1$
the fibre $Z_x$ is a single point of $C_2$; since every point of $E$ is a limit
of such points, $M_f(k)=\{k\}$ for $k\in E$. For $k\in Q'\setminus E$,
$M_f(k)\subseteq\{x\}\times[0,1]$. The function $g=\max(0,-y)$ on $S$, $g=0$ on
$Q'$, is continuous (both definitions give $0$ on $Q'\cap S=\{0\}\times[0,1]$)
and constant on every $M_f(k)$, but $g(0,-1)=1\ne0=g(0,0)$ although both points
lie in $Z_0=E$. Thus $f\notin\D(K)$ and $K\notin\DCP$.
\end{example}

\section{Individual functions in $\D(K)$}\label{s:functions}

For continua such as solenoids and the buckethandle continuum we do not
know whether they are in $\DCP$, but we can show that many individual
functions belong to $\D(K)$. The main result, Theorem~\ref{t:nonconst}, covers
every function that is constant on no arc, on a continuum all of whose arc
components are dense. Retractions onto locally connected continua give further
functions, with no assumption on $L(K)$ (Corollary~\ref{c:retract}); we apply
this to radial functions on smooth dendroids and on cones. We then collect
what these results give for solenoids, the buckethandle continuum and fans
with the property of Kelley.

\begin{theorem}\label{t:nonconst}
Let $K$ be a metrizable continuum all of whose arc components are dense, and
let $f\in C(K)$ be constant on no arc. Then $f\in\D(K)$.
\end{theorem}

\begin{proof}
We may assume that $K$ is not a point. Every nonempty open set $U$ meets a
dense, nondegenerate arc component, hence contains an arc (an initial piece of
an arc leaving a point of $U$); so all fibres have empty interior. Let
$\min f<t<\max f$. A dense arc component contains points with $f<t$ and points
with $f>t$; let $\gamma$ be an arc joining two such points. Since $f$ is not
constant on subarcs of $\gamma$, the relative closures of $\gamma\cap\{f<t\}$
and $\gamma\cap\{f>t\}$ cover $\gamma$ and meet at a point of the relative
interior of $\gamma$, which is a crossing point at level $t$; by
Lemma~\ref{l:arc} every $t\in(\min f,\max f)$ is collapsed. If $u\in Z_c^+$
then $c<\max f$, and $u$ is the limit of points $u_n$ with $c<f(u_n)<\max f$;
similarly for $Z_c^-$. Thus (a) of Theorem~\ref{t:crit} holds, (b) holds by
Lemma~\ref{l:dense-arc}, and Theorem~\ref{t:crit} applies to every level.
\end{proof}

\begin{remark}
Theorem~\ref{t:nonconst} contains~\cite[Theorem~1.19(2)]{Ki10} (arc connected
$K$ and totally disconnected fibres), whose proof was omitted there. Besides
arc connected continua, it applies to solenoids and the buckethandle
continuum; see below.
\end{remark}

\subsection*{Retractions}

Recall from Section~1 the notions of a dendroid and of a smooth dendroid.
Charatonik and Eberhart~\cite[Theorem~10]{CE} proved that $K$ is smooth with
initial point $p$ if and only if $K$ admits a metric $d$ that is \emph{radially
convex} with respect to $p$: $x\in py$, $x\ne y$ imply $d(p,x)<d(p,y)$.

\begin{proposition}\label{p:smooth}
Let $K$ be a smooth dendroid with initial point $p$, $d$ a metric on $K$ that is
radially convex with respect to $p$, and $\rho=d(p,\cdot)$. Then
$\varphi\circ\rho\in\D(K)$ for every $\varphi\in C(\rho(K))$.
\end{proposition}

\begin{proof}
Let $R=\max\rho$ and choose $e$ with $\rho(e)=R$. By radial convexity $\rho$ is
strictly increasing along the arc $pe$, so it maps $pe$ homeomorphically onto
$[0,R]=\rho(K)$; let $\gamma:[0,R]\to pe$ be the inverse. Then
$r=\gamma\circ\rho$ is a retraction of $K$ onto the arc $pe$
(cf.~\cite[Theorem~1]{Lum76}) and
$\rho\circ r=\rho$, hence $f\circ r=f$ for $f=\varphi\circ\rho$. Apply
Corollary~\ref{c:retract}.
\end{proof}

Note that $\varphi$ is arbitrary; it may be constant on intervals or behave like
the Cantor function, which is where the lemmas of Section~\ref{s:lemmas} are
most difficult to apply. (For a cone $\operatorname{Cone}(Y)$ the same
argument, with the retraction $(y,t)\mapsto(y_0,t)$ onto a single tooth, gives
every function of $t$ alone; but by Theorem~\ref{t:combfan} every cone is in
$\DCP$ anyway.)

Retractions that approximate the identity give a dense set of functions in
$\D(K)$.

\begin{proposition}\label{p:retdense}
Let $K$ be a compact metrizable space, and let $r_n$ be retractions of $K$ onto
closed subsets $L_n\in\DCP$ such that $r_n\to\mathrm{id}_K$ uniformly. Then
$\varphi\circ r_n\in\D(K)$ for every $\varphi\in C(K)$ and every $n$; in
particular, $\D(K)$ is dense in $C(K)$.
\end{proposition}

\begin{proof}
For $\varphi\in C(K)$ put $f=\varphi\circ r_n$. Then $f\circ r_n=f$ and
$L_n\in\DCP$, so $f\in\D(K)$ by Corollary~\ref{c:retract}. Since $\varphi$ is
uniformly continuous, $\varphi\circ r_n\to\varphi$ uniformly.
\end{proof}

\begin{corollary}\label{c:invret}
Let $K=\varprojlim(X_n,f_n)$, where $X_1\subseteq X_2\subseteq\cdots$ are
compact metrizable spaces in $\DCP$ and each $f_n:X_{n+1}\to X_n$ is a
retraction. Then $\D(K)$ is dense in $C(K)$.
\end{corollary}

\begin{proof}
Let $d_m$ be a metric on $X_m$ bounded by $1$ and
$d(x,y)=\sum_m2^{-m}d_m(x_m,y_m)$ on $K$. For $x\in X_n$ put
$e_n(x)_m=f_m\circ\dots\circ f_{n-1}(x)$ for $m<n$ and $e_n(x)_m=x$ for
$m\ge n$; this is a thread because $f_m$ is the identity on $X_m$. Thus
$e_n:X_n\to K$ is an embedding with $\pi_n\circ e_n=\mathrm{id}$, and
$r_n=e_n\circ\pi_n$ is a retraction of $K$ onto $e_n(X_n)\cong X_n$. Since
$r_n(y)$ and $y$ have the same first $n$ coordinates,
$d(r_n(y),y)\le2^{-n}$. Apply Proposition~\ref{p:retdense}.
\end{proof}

Such inverse limits occur, for instance, in the construction of smooth
dendroids in~\cite[proof of Theorem~9.2]{Ni89}.

By a theorem of Fugate~\cite[Theorem~2]{Fu71b}, for every smooth dendroid $X$
and every $\varepsilon>0$ there are a tree (a finite dendrite) $T\subseteq X$
and a retraction $r:X\to T$ that moves no point of $X$ as much as
$\varepsilon$. Trees are locally connected, hence in $\DCP$
(Corollary~\ref{c:old}), so Proposition~\ref{p:retdense} gives the following.

\begin{corollary}\label{c:smoothdendroid}
Let $X$ be a smooth dendroid. Then $\D(X)$ is dense in $C(X)$.
\end{corollary}

For fans smoothness is not needed: by~\cite[Theorem~1]{Fu71a}, for every fan
$X$ and every $\varepsilon>0$ there are a finite fan $Y\subseteq X$ and a
retraction $r:X\to Y$ with $d(x,r(x))<\varepsilon$ for all $x\in X$. Finite fans
are locally connected, so Proposition~\ref{p:retdense} again applies.

\begin{corollary}\label{c:allfans}
Let $X$ be a fan. Then $\D(X)$ is dense in $C(X)$.
\end{corollary}

\begin{remark}\label{r:density}
Density of $\D(K)$ says nothing about whether $K\in\DCP$. The $\sin(1/x)$
continuum $S$ is not in $\DCP$ (Example~\ref{e:sin}), yet it admits retractions onto subarcs that
move no point by more than $1/n$ (see the introduction of~\cite{EF71}): choose
an interval $[a,a']\subseteq(0,1/n]$ on which $\sin(1/x)$ increases from $-1$
to $1$, send every point of $S$ with first coordinate at most $a$ to the point
of the graph over $[a,a']$ with the same second coordinate, and fix the rest of
the graph. This is a retraction of $S$ onto the graph over $[a,1]$, an arc.
Arcs are in $\DCP$, so $\D(S)$ is dense in $C(S)$ by
Proposition~\ref{p:retdense}.
\end{remark}

For smooth fans the retractions can be written down directly, and they can be
chosen compatibly (Lemma~\ref{l:fanlim} below). A fan is smooth if and only if it
embeds in the Cantor fan~\cite[Theorem~9]{Ch67}, \cite[Corollary~4]{Eb69}
(cf.\ the proof of~\cite[Theorem~6.9.3]{Ma05}).

\begin{corollary}\label{p:smoothfan}
Let $X$ be a smooth fan. Then there are retractions $r_n$ of $X$ onto finite
unions of legs of $X$ such that $r_n\to\mathrm{id}_X$ uniformly, and
$\varphi\circ r_n\in\D(X)$ for every $\varphi\in C(X)$ and every $n$. In
particular, $\D(X)$ is dense in $C(X)$.
\end{corollary}

\begin{proof}
Regard $X$ as a subcontinuum of $F=\operatorname{Cone}(C)$ and write $[x,t]$ for
the image of $(x,t)\in C\times[0,1]$. Since $X$ is not an arc, its ramification
point is the vertex $v$ of $F$, the only ramification point of $F$; hence for
every $y\in X$ the arc $vy$ is an initial segment of a leg of $F$. Thus
$X=\{[x,t]:t\le\ell(x)\}$, where $\ell(x)=\max\{t:[x,t]\in X\}$; the function
$\ell:C\to[0,1]$ is upper semicontinuous because $X$ is closed.

Let $\mathcal P_n$ be a partition of $C$ into finitely many clopen sets of
diameter less than $1/n$. In each $P\in\mathcal P_n$ choose $x_P$ with
$\ell(x_P)=\max_P\ell$ (the maximum exists since $\ell$ is upper semicontinuous
and $P$ is compact), and define $r_n[x,t]=[x_P,t]$ for $x\in P$. This is
well defined at $v$, continuous because $x\mapsto x_P$ is locally constant, and
maps $X$ into $X$ because $t\le\ell(x)\le\ell(x_P)$. It is the identity on the
legs $\{[x_P,t]:t\le\ell(x_P)\}$, $P\in\mathcal P_n$, whose union $L_n$ is a
finite union of arcs with a common end point, hence a locally connected
continuum, so $L_n\in\DCP$ by Corollary~\ref{c:old}. By uniform continuity of
the quotient map $C\times[0,1]\to F$, $r_n\to\mathrm{id}_X$ uniformly. Apply
Proposition~\ref{p:retdense}.
\end{proof}

The retractions of Corollary~\ref{p:smoothfan} can be chosen compatibly, which
represents every smooth fan as an inverse limit of the kind considered in
Corollary~\ref{c:invret}.

\begin{lemma}\label{l:fanlim}
Let $X$ be a smooth fan. Then $X$ is homeomorphic to $\varprojlim(L_n,s_n)$,
where each $L_n$ is a finite union of arcs with a common end point, $L_n
\subseteq L_{n+1}$, and $s_n:L_{n+1}\to L_n$ is a retraction.
\end{lemma}

\begin{proof}
We use the notation of the proof of Corollary~\ref{p:smoothfan}, choosing the
partitions nested: every member of $\mathcal P_{n+1}$ is contained in a member
of $\mathcal P_n$. We choose the points $x_P$ inductively. If $P\in\mathcal P_n$
and $P'\in\mathcal P_{n+1}$ is the member with $x_P\in P'\subseteq P$, put
$x_{P'}=x_P$; this is admissible since $\ell(x_P)=\max_P\ell\ge\max_{P'}\ell$.
For the other members of $\mathcal P_{n+1}$ choose $x_{P'}$ as before.

Then $L_n\subseteq L_{n+1}$. Moreover $r_n\circ r_{n+1}=r_n$: if $x\in
P'\subseteq P$ with $P'\in\mathcal P_{n+1}$, $P\in\mathcal P_n$, then $x_{P'}\in
P$, so $r_n(r_{n+1}[x,t])=r_n[x_{P'},t]=[x_P,t]=r_n[x,t]$. Hence
$s_n=r_n|_{L_{n+1}}$ is a retraction of $L_{n+1}$ onto $L_n$, and by induction
$r_n\circ r_m=r_n$ for $m\ge n$.

Define $h:X\to\varprojlim(L_n,s_n)$ by $h(y)=(r_n(y))_n$; this is a thread
because $r_n\circ r_{n+1}=r_n$. The map $h$ is continuous, and injective because
$r_n(y)\to y$. Let $(z_n)$ be a thread and $m\ge n$. Since
$r_n=r_n\circ r_{n+1}\circ\dots\circ r_{m-1}$ and this composition agrees with
$s_n\circ\dots\circ s_{m-1}$ on $L_m$, we get $r_n(z_m)=z_n$. If $y$ is a
limit point of
$(z_m)$ in $X$, continuity of $r_n$ gives $r_n(y)=z_n$ for every $n$, i.e.\
$h(y)=(z_n)$. Thus $h$ is a continuous bijection of compact metrizable spaces,
hence a homeomorphism.
\end{proof}

Corollary~\ref{p:smoothfan} applies to the Lelek fan, to which
Corollary~\ref{c:fans} below does not apply (see Problem~\ref{p:dendroid}).
As Remark~\ref{r:density} shows, density of $\D(X)$ gives no evidence that
$X\in\DCP$.

\begin{example}\label{e:conv}
The converse of Lemma~\ref{l:retract} fails. Let $S$ and $E$ be as in
Example~\ref{e:sin}, $A=[-1,0]\times\{0\}$, $K=S\cup A$, and let $r:K\to S$ be
the identity on $S$ and map $A$ to $(0,0)$. Put $f(x,y)=x$ on $S$ and
$f=f|_S\circ r$, so $f=0$ on $A$ and $Z_0=E\cup A$. By Example~\ref{e:sin},
$g_S(x,y)=y$ is $f|_S$-admissible. In $K$, the set $A\setminus\{(0,0)\}$ is open
and contained in $Z_0$, so $M_f(k)=Z_0$ for $k\in A\setminus\{(0,0)\}$ by
Lemma~\ref{l:interior}. Hence $g=g_S\circ r$ is not $f$-admissible, although
$g=g\circ r$ and $g|_S=g_S$. Attaching an arc to a fibre enlarges $M_f$.
\end{example}

\subsection*{Solenoids and the buckethandle continuum}

Recall that a \emph{solenoid} is the inverse limit of circles with bonding
maps $z\mapsto z^{n_i}$, $n_i\ge2$ (see~\cite[2.8]{Na92} for the case
$n_i=p$), and the \emph{buckethandle continuum} is the inverse limit of copies
of $[0,1]$ with the tent map $t\mapsto1-|2t-1|$~\cite[2.9]{Na92}. The tent
map is open, and by~\cite[Theorem~232]{IM12} the same continuum is obtained
from any sequence of maps of $[0,1]$ onto itself each of which is the union of
two homeomorphisms of complementary subintervals onto $[0,1]$. (The
buckethandle is often called the Knaster continuum; it should not be confused
with Knaster's hereditarily indecomposable continuum of 1922, which is the
pseudo-arc~\cite[p.~229]{Na92} and contains no arcs.) Both are
indecomposable, and in both the arc components coincide with the composants:
for solenoids because every nondegenerate proper subcontinuum is an arc
(\cite[2.16]{Na92} for $n_i=p$; the general case is the same), and for the
buckethandle see~\cite[p.~268]{St02}.
Since composants are dense~\cite[5.20(a)]{Na92}, all arc components are
dense.

Hence, by Theorem~\ref{t:nonconst}, on a solenoid and on the buckethandle
continuum every function that is constant on no arc belongs to $\D(K)$; for
solenoids see also Proposition~\ref{p:sol} below.

\begin{remark}\label{r:confluent}
The density of arc components is not special to these examples. Call a
continuum $K$ \emph{confluently locally connected-like} if for every
$\varepsilon>0$ there is a confluent $\varepsilon$-mapping of $K$ onto a
locally connected continuum; this holds, for instance, for every inverse limit
of locally connected continua with confluent (e.g.\ open, or monotone) bonding
maps. By~\cite[Corollary~2.26]{CCP08}, every such continuum is arc Kelley, and
in particular all its arc components are dense. Hence, on every such
continuum, every function that is constant on no arc belongs to $\D(K)$
(Theorem~\ref{t:nonconst}), and condition (b) of Theorem~\ref{t:crit} holds
for every $f$ and every level with empty interior (Lemma~\ref{l:dense-arc});
since levels with interior are collapsed (Lemma~\ref{l:interior}), whether
$K\in\DCP$
depends only on condition (a), as for arc connected continua. By a theorem of
Oversteegen and Prajs~\cite[Corollary~3.5]{OP03} (see also~\cite[Section~3]{CCP08}), the nondegenerate atriodic
continua of this kind are exactly the solenoids and the Knaster type continua.
\end{remark}

\begin{remark}\label{r:indec}
Indecomposability alone does not place a continuum in $\DCP$. By a theorem of
Fukaishi and Matsuhashi~\cite[Theorem~1.1]{FM17}, for every continuum $X$ there are an
indecomposable continuum $Y\supseteq X$ and an open retraction $r:Y\to X$
(whose fibres are Cantor sets). If $Y$ were in $\DCP$, so would be its open
image $X$ (Theorem~\ref{t:open-map}). Taking for $X$ the $\sin(1/x)$ continuum
(Example~\ref{e:sin}), we obtain an indecomposable continuum that is not in
$\DCP$; more generally, every continuum not in $\DCP$ is an open retract of an
indecomposable continuum not in $\DCP$. Since the $\sin(1/x)$ continuum is
arc-like, by~\cite[Theorem~3.3(2)]{MOB25} the continuum $Y$ can moreover be
chosen tree-like, so there are tree-like indecomposable continua outside
$\DCP$.
\end{remark}

\begin{proposition}\label{p:sol}
Let $K$ be a solenoid, let $f\in C(K)$, let $I$ be the union of the interiors
of the fibres of $f$, and let $\mathrm{TF}$ be the set of points having a
Cantor box neighbourhood on each of whose arcs $\{x\}\times J$ the function
$f$ is constant. If $\mathrm{TF}\subseteq\cl I$, then $f\in\D(K)$.
\end{proposition}

\begin{proof}
Levels with interior are collapsed. Let $\Int Z_c=\emptyset$; (b) holds by
Lemma~\ref{l:dense-arc}. For (a) let $u\in Z_c^+$ and $U$ a neighbourhood of
$u$. If $U\cap\{f>c\}$ meets $I$ we get a point with $f>c$ on a collapsed
level. Otherwise $W'=U\cap\{f>c\}$ is a nonempty open set covered by Cantor
boxes (every point of a solenoid has Cantor box neighbourhoods) and disjoint
from $\cl I\supseteq\mathrm{TF}$, so $f$ is nowhere locally constant on $W'$
and some box in $W'$ has an arc on which $f$ is not constant. That arc carries
crossing points (at all but countably many levels of its image), which are
collapsed by Proposition~\ref{p:box}. Apply Theorem~\ref{t:crit}.
\end{proof}

\subsection*{Fans with the property of Kelley}

A continuum $X$ with metric $d$ has the \emph{property of Kelley} if for every
$\varepsilon>0$ there is $\delta>0$ such that whenever $a,b\in X$,
$d(a,b)<\delta$, and $A$ is a subcontinuum of $X$ containing $a$, there is a
subcontinuum $B$ of $X$ containing $b$ with Hausdorff distance
$H(A,B)<\varepsilon$; see~\cite{CC88b}. An inverse limit of continua with
the property of Kelley with confluent (in particular, open) bonding maps has
the property of Kelley~\cite[Theorem~219]{IM12}; for fans there is a converse
(Corollary~\ref{c:fans} below).
Corollary~\ref{c:smooth} covers the
smooth fans with a closed set of end points. For fans with the property of
Kelley, Theorem~\ref{t:open-map} gives the following.

\begin{corollary}\label{c:fans}
Let $X$ be a fan with the property of Kelley. Then
$X=\varprojlim(X_n,f_n)$ with $X_n$ locally connected fans and $f_n$ open
\cite[Theorem~3]{CC88b}, and with $\pi_n:X\to X_n$ the projections, the
functions $h\circ\pi_n$, $h\in C(X_n)$, $n\in\N$, belong to $\D(X)$ and are
dense in $C(X)$.
\end{corollary}

\begin{proof}
The projections $\pi_n$ are surjective~\cite[Exercise~2.14]{Na92} and open
(\cite[Corollary~2.1.10]{Ma05}; we recall the argument):
the sets $\pi_m^{-1}(U)$, $m\ge n$, $U$ open in $X_m$, form a base of $X$, and
$\pi_n(\pi_m^{-1}(U))=f_{n,m}(U)$ is open, $f_{n,m}=f_n\circ\dots\circ
f_{m-1}$ being open. Since $X_n$ is a locally connected continuum, $h\in
\D(X_n)$ by Corollary~\ref{c:old}, so $h\circ\pi_n\in\D(X)$ by
Theorem~\ref{t:open-map}. Since $\pi_n=f_n\circ\pi_{n+1}$, these functions
form an algebra; it contains the constants and separates the points of $X$, so
it is dense in $C(X)$ by the Stone--Weierstrass theorem.
\end{proof}

\section{Open problems}\label{s:problems}

\begin{problem}\label{p:arc}
Is every arc connected metrizable continuum in $\DCP$? By the remark preceding
Theorem~\ref{t:arc-dense}, condition (b) of Theorem~\ref{t:crit} always holds
for arc connected $K$, so only condition (a) has to be verified. By
Theorem~\ref{t:nonconst} the difficulty comes only from functions that are
constant on some arc, and by Lemma~\ref{l:interior} only from their levels
with empty interior.

A weaker question: is every arc connected metrizable continuum in which the
points having a locally connected open neighbourhood or a \emph{Cantor} box
neighbourhood are dense in $\DCP$? This would contain
Theorems~\ref{t:arc-dense} and~\ref{t:grid}, and Theorem~\ref{t:combfan} for
the comb and the Cantor fan. The restriction to Cantor boxes matters: with
boxes as in Remark~\ref{r:box} the question would be equivalent to the
problem, since for an arc connected continuum $Y$ every point of $Y\times
(0,1)$ has a box neighbourhood, so a positive answer would give
$Y\times[0,1]\in\DCP$ and hence $Y\in\DCP$ by Corollary~\ref{c:product}. For
continua with tooth charts, Theorem~\ref{t:teeth} reduces the question to
showing that the roots lie in $\cl\Sigma$. A test case is the dendroid
of~\cite[Example~7.15]{Il25}: $C\times[0,1]$ with the two end points of each
gap of $C$ in $C\times\{0\}$ identified, whose roots form an arc.
\end{problem}

\begin{problem}\label{p:dendroid}
Is every dendroid in $\DCP$? Every smooth dendroid? Every fan? No dendroid
outside $\DCP$ is known. By Theorem~\ref{t:combfan} and
Corollary~\ref{c:smooth}, the Cantor comb, the Cantor fan and all smooth fans
with a closed set of end points are in $\DCP$. A natural test case is the
Lelek fan, the unique smooth fan with a dense set of end
points~\cite{Ch89,BO90}. Since its set of end points is dense, it does not
have the property of Kelley~\cite[Theorem~6]{CC89}; in particular it is
neither a retract nor an open image of the Cantor fan~\cite[Theorems~6,
8]{CC89}, and Corollary~\ref{c:fans} does not apply to it. On the other hand,
every Cantor fan contained in the Lelek fan is a retract of
it~\cite[Theorem~3.1]{BEGK}, so by Corollary~\ref{c:retract} and
Theorem~\ref{t:combfan} every function on the Lelek fan that factors through
such a retraction belongs to $\D(X)$, where $X$ is the Lelek fan, and by
Corollary~\ref{p:smoothfan} $\D(X)$ is dense in $C(X)$; by
Corollaries~\ref{c:smoothdendroid} and~\ref{c:allfans} the same holds for
every smooth dendroid and every fan; by Remark~\ref{r:density}, however, this
is no evidence for $\DCP$. Since its end points are dense and no end point of
a fan lies in a Cantor box (every point of a box is interior to an arc of the
box), it contains no Cantor boxes, and Proposition~\ref{p:box} does not apply.
A convenient explicit model for constructing test functions is given
in~\cite[Theorem~4.34]{BEK23}: if $0<r<1<\rho$ and $r^k=\rho^\ell$ only for
$k=\ell=0$, then the set of all sequences $(x,a_1x,a_1a_2x,a_1a_2a_3x,\dots)
\in[0,1]^{\N}$ with $a_i\in\{r,\rho\}$ is a Lelek fan, the union of straight
segments from the origin of the Hilbert cube; it is the generalized inverse
limit of $[0,1]$ with the single set-valued bonding function whose graph is the
union of the segments $y=rx$ and $y=\rho x$ in $[0,1]^2$, and it is mapped into
itself by the shift. Every coordinate function is strictly increasing along
each leg and hence in $\D(X)$ by Theorem~\ref{t:nonconst}; a function outside
$\D(X)$, if one exists, must be constant on some arcs. The Lelek fan also
arises in complex dynamics: for $0<\lambda<1/e$ the Julia set $J_\lambda$ of
$z\mapsto\lambda e^z$, together with the point $\infty$, is homeomorphic to the
Lelek fan, with $\infty$ as its top~\cite{AO93} (see also~\cite{BO90}). Thus the
question is whether $\{M_f\}''=\A_f$ for every multiplication operator $M_f$ on
$C(J_\lambda\cup\{\infty\})$.

By~\cite{Ch84}, a dendroid is smooth if and only if it is the inverse limit of
finite dendrites with bonding maps that are monotone relative to the points
of a thread. These bonding maps cannot be chosen open, nor even confluent, in
general: locally connected continua have the property of Kelley
(see~\cite[p.~1105]{CCP08}), so by~\cite[Theorem~219]{IM12} an inverse limit
of dendrites with confluent (for instance open, or monotone in the usual sense)
bonding maps has the property of Kelley, which the Lelek fan does not.
Thus the first question of Problem~\ref{p:invlim} does not reach the Lelek
fan (the second one does, by Lemma~\ref{l:fanlim}), and it is natural to
ask whether $\DCP$ is preserved under inverse limits with bonding maps that
are monotone relative to the points of a thread.
\end{problem}

\begin{problem}
Is a retract of a continuum in $\DCP$ in $\DCP$? Is a continuous image of a
continuum in $\DCP$ in $\DCP$? By~\cite[Theorem~6]{CC89}, every fan with the
property of Kelley is a retract of the Cantor fan, so by
Theorem~\ref{t:combfan} a positive answer to the first question would put all
of them in $\DCP$. By Kuperberg's theorem~\cite[Theorem~3.5]{Ku} (see
also~\cite[Theorem~12]{CC89}) the continuous images of the Cantor fan are
exactly the uniformly pathwise connected continua; a positive answer to the
second question would put all of them in $\DCP$.
Example~\ref{e:conv} shows that Lemma~\ref{l:retract} alone does not settle
the first question. By Theorem~\ref{t:open-map} both questions have a positive
answer for open retractions and open maps.
\end{problem}

\begin{problem}\label{p:invlim}
Is $\DCP$ preserved under inverse limits with open bonding maps? A positive
answer would give, via Corollary~\ref{c:fans}, that all fans with
the property of Kelley are in $\DCP$; since the
bonding maps in the definitions of solenoids and of the buckethandle continuum
(Section~\ref{s:functions}) are open, it would also put all solenoids
and the buckethandle continuum in $\DCP$. The same question can be asked for
inverse limits whose bonding maps are retractions, as in
Corollary~\ref{c:invret}; for these $\D(K)$ is always dense in $C(K)$. By
Lemma~\ref{l:fanlim}, every smooth fan, in particular the Lelek fan, is such a
limit of finite unions of arcs with a common end point.
\end{problem}

\begin{problem}\label{p:converse}
Find a characterization of $\DCP$ in terms of the topology of $K$ alone, that
is, one that does not require verifying $f\in\D(K)$ for every $f\in C(K)$
separately.
In particular, is the converse of Theorem~\ref{t:fat} true: if $L(K)$ is dense
and $K\in\DCP$, is every arc component of $K$ either dense or contained in the
closure of its interior? This holds in all examples we know: the $\sin(1/x)$
continuum, Examples~\ref{e:sin-arc} and~\ref{e:halfsquare}, the union of
$\{0\}\times[-1,1]$ with two $\sin(1/x)$ curves approaching it from both sides,
and a ray spiralling onto a circle (the first and the last two are not in
$\DCP$ by Proposition~\ref{p:ray}). A positive answer would give a purely
topological characterization of $\DCP$ for continua with $L(K)$ dense.
\end{problem}

\subsection*{Acknowledgment}
Several results of this paper, among them the topological reformulation in
Section~\ref{s:char}, the sufficient condition for collapsed levels in Section~\ref{s:main} (Theorem~\ref{t:crit}) and the theorems on
the Cantor comb, the Cantor fan and the Cantor grid, were obtained in an extended interaction between the author
and Claude, an AI system developed by Anthropic. Claude also prepared drafts of
the text. A separate instance of Claude, given only the manuscript, reviewed the
proofs. The author has checked all arguments and is fully responsible for the
contents of the paper.

\end{document}